\documentclass[review]{elsarticle}
\usepackage[T1]{fontenc}
\usepackage[utf8]{inputenc}
\usepackage{units}
\usepackage{mathtools}
\usepackage{amsmath}
\usepackage{amssymb}
\usepackage{stmaryrd}
\usepackage{graphicx}
\usepackage[english]{babel}
\usepackage{hyperref}

\makeatletter

\providecommand{\tabularnewline}{\\}

\numberwithin{equation}{section}
\newtheorem{rem}{\protect\remarkname}
\newtheorem{defn}{\protect\definitionname}

\newtheorem{thm}{Theorem}
\newtheorem{lem}[thm]{Lemma} 
\newdefinition{rmk}{Remark} 
\newproof{pf}{Proof}
\newproof{pot}{Proof of Theorem \ref{thm2}}

\renewcommand{\vec}[1]{\mbox{\boldmath$#1$}}

\usepackage{xfrac} 
\usepackage{algorithmic}

\renewcommand\[{\begin{equation}}
\renewcommand\]{\end{equation}}

\makeatother

\providecommand{\definitionname}{Definition}

\providecommand{\remarkname}{Remark}

\begin{document}

\begin{frontmatter}

\title{A Hybrid Method and Unified Analysis of Generalized \\
Finite Differences and Lagrange Finite Elements}

\author[1]{Rebecca Conley}
\ead{rconley@saintpeters.edu}

\author[2,3]{Tristan J. Delaney}
\ead{tristan.delaney@synopsys.com}

\author[2]{Xiangmin Jiao\corref{cor1}}
\ead{xiangmin.jiao@stonybrook.edu}

\cortext[cor1]{Corresponding author}
\address[1]{Dept. of  Math., Saint Peter's University, Jersey City, NJ 07306, USA.}
\address[2]{Dept. of Applied Math. \& Stat. and Institute for Advanced Computational Science, Stony Brook University, Stony Brook, NY 11794, USA.}
\address[3]{Current address: Synopsys, Inc., Mountain View, CA 94043, USA.}

\begin{abstract}
Finite differences, finite elements, and their generalizations are
widely used for solving partial differential equations, and their
high-order variants have respective advantages and disadvantages.
Traditionally, these methods are treated as different (strong vs.
weak) formulations and are analyzed using different techniques (Fourier
analysis or Green's functions vs. functional analysis), except for
some special cases on regular grids. Recently, the authors introduced
a hybrid method, called \emph{Adaptive Extended Stencil FEM} or \emph{AES-FEM}
(\emph{Int. J. Num. Meth. Engrg., 2016}, DOI:10.1002/nme.5246), which
combines features of \emph{generalized finite differences} and \emph{Lagrange
finite elements} to achieve second-order accuracy over unstructured
meshes. However, its analysis was incomplete due to the lack of existing
mathematical theory that unifies the formulations and analysis of
these different methods. In this work, we introduce the framework
of \emph{generalized weighted residuals} to unify the formulation
of finite differences, finite elements, and AES-FEM. In addition,
we propose a unified analysis of the \emph{well-posedness}, \emph{convergence},
and \emph{mesh-quality dependency} of these different methods. We
also report numerical results with AES-FEM to verify our analysis.
We show that AES-FEM improves the accuracy of generalized finite differences
while reducing the mesh-quality dependency and simplifying the implementation
of high-order finite elements.
\end{abstract}
\begin{keyword}
partial differential equations\sep finite element methods\sep generalized
finite differences\sep generalized weighted residuals\sep stability\sep convergence
\MSC[2010] 65N06 \sep  65N30 \sep 65N12
\end{keyword}
\end{frontmatter}


\section{Introduction}

Finite differences, finite elements, and their generalizations are
widely used for solving partial differential equations (PDEs), or
more precisely, for the spatial discretization for their associated
boundary value problems (BVPs). The finite difference methods (FDM)
are standard techniques in numerical analysis \cite{leveque2007finite,strikwerda2004finite}
and are widely used for solving hyperbolic PDEs in computational fluid
dynamics. The finite element methods (FEM), on the other hand, are
the most successful method for solving elliptic and parabolic PDEs
(e.g., \cite{strang1973analysis,zienkiewicz2013finite}). Since many
PDEs, such as advection-diffusion-reaction equations, Navier-Stokes
equations, etc., are multiphysics in nature, involving both parabolic
(elliptic) and hyperbolic components, it has been of great interest
for applied mathematicians to develop hybrid methods that combine
the advantages of FEM and FDM. The most notable examples are \emph{discontinuous
Galerkin} methods \cite{Cockburn2000,Riviere2008DG} and some \emph{finite
volume methods} \cite{leveque2002finite,li2000generalized}, which
use discontinuous test functions analogous to the \emph{nonconforming
finite elements }\cite[Section 10.3]{brenner2008mathematical} and
use finite-difference style computations of interface fluxes or jump
conditions (such as WENO \cite{shu2003high,costa2007high} and Lax–Friedrichs
limiters \cite[p. 199]{leveque2002finite}) along element or cell
boundaries.

In \cite{conley_2016,conley2016overcoming}, the authors introduced
a new hybrid method, called the \emph{adaptive extended stencil finite
element method} or \emph{AES-FEM}. Like Lagrange finite elements,
AES-FEM has $C^{0}$ continuous test functions, so there are no explicit
interface and jump conditions in its\emph{ }variational forms, unlike
nonconforming finite elements. However, unlike finite elements, AES-FEM
uses least-squares based trial functions similar to those of generalized
finite differences \cite{Benito2007}. We refer to these trial functions
as \emph{generalized Lagrange polynomial} (GLP) basis functions, which
are not \emph{$C^{0}$ continuous} but have similar properties to
Lagrange interpolation. The GLP basis functions introduce a ``variational
crime'' (cf. \ref{subsec:Variational-Crimes}), which is similar
to, but yet different from, that of other nonconforming finite elements.
The analysis of GLP-based methods, including generalized finite difference
method (GFDM) and AES-FEM, requires a mathematical analysis that unifies
the classical analysis of finite differences and finite elements.
This unification is the primary goal of this work, which will reveal
some insights from a theoretical point of view, and also enable a
rigorous generalization of AES-FEM to higher-order accuracy from a
practical point of view.

The main contributions of this work are as follows. First, we unify
the formulations of GFDM, FEM, and AES-FEM under the framework of
\emph{generalized weighted residuals} (GWR), of which the trial functions
are either Lagrange or generalized Lagrange basis functions. Second,
we establish the conditions for \emph{well-posedness}, \emph{convergence},\emph{
}and\emph{ superconvergence }of GFDM and AES-FEM, and compare their
mesh-quality requirements against Lagrange finite elements. Third,
we prove and also demonstrate the high-order convergence of AES-FEM.
For simplicity, we assume exact geometry for Neumann boundaries in
this paper, and we defer the treatment of Neumann boundary conditions
over approximate curved boundaries to future work. 

The remainder of the paper is organized as follows. Section~\ref{sec:Background}
reviews the (G)FDM and FEM for boundary value problems, as well as
their respective classical analyses. Section~\ref{sec:Generalized-Weighted-Residuals}
introduces the concept of \emph{generalized Lagrange polynomial} basis
functions and the framework of \emph{generalized weighted residuals}
(GWR), which unifies the formulations of GFDM, FEM, and AES-FEM. Section~\ref{sec:well-posedness}
analyzes the well-posedness of GWR methods. Section~\ref{sec:convergence}
addresses the convergence of GFDM and AES-FEM to confirm our analysis.
Section~\ref{sec:numerical-results} presents some numerical results.
Section~\ref{sec:Conclusions} concludes this paper with a discussion
on future work.

\section{\label{sec:Background} Background and Preliminaries}

In this section, we briefly review finite differences, finite elements,
and some of their generalizations for boundary value problems. We
refer readers to \cite{thomee2001finite} for a brief history of these
different methods. For completeness, we review some relevant details
about the generalized Lagrange polynomials, AES-FEM, and functional
analysis in the appendices.

\subsection{Boundary Value Problems}

Let $\Omega\subset\mathbb{R}^{d}$ be a bounded, piecewise smooth
domain with boundary $\Gamma=\Gamma_{D}\cup\Gamma_{N}$, where $d$
is typically 2 or 3, and $\Gamma_{D}$ and $\Gamma_{N}$ denote the
Dirichlet and Neumann boundaries, respectively. Let $\mathcal{L}$
be a second-order linear differential operator. In general, $\mathcal{L}$
has the form of
\begin{equation}
\mathcal{L}u=-\vec{\nabla}\cdot\left(\mu\vec{\nabla}u\right)+\vec{\nu}\cdot\vec{\nabla}u+\omega^{2}u,\label{eq:2nd-order-elliptic-operator}
\end{equation}
where $\vec{\nabla}$ and $\vec{\nabla}\cdot$ denote the gradient
and divergence operators, $\mu(\vec{x}):\Omega\rightarrow\mathbb{R^{+}}$
corresponds to a diffusion coefficient, $\vec{\nu}:\Omega\rightarrow\mathbb{R}^{d}$
is a velocity field, and $\omega:\Omega\rightarrow\mathbb{R}$ is
a wavenumber or frequency. Typically, $\vec{\nabla}\cdot\vec{\nu}=0$.
A second-order partial differential equation has the form of
\begin{align}
\mathcal{L}u & =f\qquad\text{on }\Omega,\label{eq:bndy-value-problem}
\end{align}
where $f:\Omega\rightarrow\mathbb{R}$ is a source term. This general
form is known as the \emph{advection-diffusion-reaction} equations
for vector-valued PDEs.  For simplicity, we focus on scalar fields
and assume \emph{diffusion dominance} (i.e., $\mu(\vec{x})\geq C\left\Vert \vec{\nu}(\vec{x})\right\Vert h$
for some $C\ge1$, where $h$ denotes some characteristic edge length
of the mesh; cf. \ref{subsec:Functional-Analysis}). If $\omega=0$,
then its corresponding PDE is the \emph{advection-diffusion equation}.
A\emph{ boundary value problem} (BVP) corresponding to the above PDE
may have some \emph{Dirichlet} \emph{boundary conditions}
\begin{align}
u & =u_{D}\qquad\text{on }\Gamma_{D},\label{eq:bndy-problem-Dirichlet}
\end{align}
and potentially some \emph{Neumann boundary conditions}
\begin{align}
\mu\partial_{\vec{n}}u & =g_{N}\qquad\text{on }\Gamma_{N},\label{eq:bndy-problem-Neumann}
\end{align}
where $\partial_{\vec{n}}$ denotes the normal derivative, i.e., $\partial_{\vec{n}}\equiv\vec{n}\cdot\vec{\nabla}$.
The boundary condition is said to be \emph{homogeneous} if $\Gamma_{D}=\Gamma$
and $u_{D}=0$.  

\subsection{\label{subsec:Finite-Differences-and} Finite Differences and Their
Generalizations}

The finite difference methods (FDM) are arguably the simplest and
the best known numerical methods for solving initial and boundary
value problems; see textbooks such as \cite{leveque2007finite,strikwerda2004finite}.
In a nutshell, FDM approximates the partial derivatives in (\ref{eq:2nd-order-elliptic-operator})
using finite difference operators, which result in a system of algebraic
equations $\vec{A}\vec{x}=\vec{b}$. The analysis of FDM is also conceptually
simple: If the truncation errors in each algebraic equation are \emph{consistent}
(i.e., $\mathcal{O}(h^{k})$ for some $k>0$ as the edge length $h$
tends to $0$) and the \emph{absolute condition number} of the algebraic
system (i.e., $\left\Vert \vec{A}^{-1}\right\Vert $) is bounded independently
of $h$ (i.e., $\mathcal{O}(1)$), then the finite difference method
converges in exact arithmetic. This condition is known as the \emph{fundamental
theorem of numerical analysis} and is simply stated as ``\emph{consistency
and stability imply convergence}'' \cite[p. 124]{strang1973analysis}.
The main argument then involves bounding the absolute condition number,
which traditionally is done using Fourier analysis on regular grids
\cite[p. 20]{leveque2007finite} or using Green's functions in 1D
\cite[p. 22]{leveque2007finite}. The effect of rounding errors is
typically omitted in the analysis of BVPs, but it has been considered
by some authors; see e.g. \cite{chan1988effectively,christiansen1994effective,li2007effective}.

The classical FDM is limited to structured meshes because the finite
difference operators are defined based on 1D polynomial interpolations
locally at each node in a dimension-by-dimension fashion. The same
approach can be utilized on curvilinear meshes for curved but relatively
simple geometries; see e.g. \cite{visbal2002use}. Some authors have
considered its generalizations to unstructured meshes or point clouds;
see e.g. \cite{Benito2007,cueto2007finite,liu2016wls,macneal1953asymmetrical}.
 Finite difference operators have been generalized to use least-squares
approximations; see e.g. \cite{Benito2007,cueto2007finite,liu2016wls}.
In this work, we use \emph{generalized finite differences} (\emph{GFD})
to refer to the least-squares-based finite difference operators and
use \emph{generalized finite difference methods} (\emph{GFDM}) to
refer to the methods that use these GFD operators to convert (\ref{eq:bndy-value-problem})
directly into algebraic equations. To the best of our knowledge, there
was no prior complete convergence analysis of GFDM for BVPs, except
for local consistency using Taylor series and the temporal aspect
of stability for time-dependent PDEs \cite{Benito2007,prieto2011application,gavete2017solving,urena2019solving}.

\subsection{\label{subsec:Finite-Elements-and}Finite Elements and Weighted Residuals}

The finite element methods (FEM) are among the most powerful and successful
methods for solving BVPs. Mathematically, FEM can be expressed using
the framework of \emph{weighted residuals} \cite{finlayson1973method},
also known as \emph{variational formulations} \cite[p. 2]{brenner2008mathematical}.
Let $\Omega_{h}$ denote the approximation of the domain $\Omega$
with a mesh. Without loss of generality, let us assume triangular
or tetrahedral meshes, and let $n$ denote the number of nodes in
$\Omega_{h}\backslash\Gamma_{D}$. Let $u_{h}$ denote the approximate
solution of the PDE on $\Omega$. The \emph{residual} of (\ref{eq:bndy-value-problem})
corresponding to $u_{h}$ is $\mathcal{L}u_{h}-f$. Let $\{\psi_{i}\mid1\leq i\le n\}$
denote the set of\emph{ test }(\emph{or weight})\emph{ functions},
which span the \emph{test space} $\Psi$. A \emph{weighted residual
method} requires the residual to be orthogonal to $\Psi$ over $\Omega$,
or equivalently, 
\begin{equation}
\int_{\Omega}\mathcal{L}u_{h}\psi_{i}\ d\vec{x}=\int_{\Omega}f\psi_{i}\ d\vec{x},\qquad\text{for }i=1,\ldots,n.\label{eq:mwr-eqn}
\end{equation}
In FEM, the test functions $\psi_{i}$ are (weakly) differentiable
and have local support.

To discretize the problem fully, let $\{\phi_{j}\mid1\leq j\leq n\}$
denote a set of basis functions, which span the \emph{trial space}
$\Phi$. Let $\vec{\Phi}$ denote the vector containing $\phi_{j}$.
We find the approximate solution $u_{h}$ in $\Phi$, i.e.,
\begin{equation}
u\approx u_{h}=\vec{\Phi}^{T}\vec{u}_{h},\label{eq:approximate-solution}
\end{equation}
where $\vec{u}_{h}$ is the unknown vector. The basis functions are
\emph{Lagrange} if $\phi_{j}(\vec{x}_{i})=\delta_{ij}$, the Kronecker
delta function; i.e., $\phi_{j}(\vec{x}_{i})=1$ if $i=j$ and $\phi_{j}(\vec{x}_{i})=0$
if $i\neq j$. With Lagrange basis functions, let $\vec{u}_{I}$ denote
the vector composed of $u(\vec{x}_{j})$. Then, $u_{I}=\vec{\Phi}^{T}\vec{u}_{I}$
is the \emph{interpolation} of $u$ in $\Phi$. Furthermore, the unknown
vector $\vec{u}_{h}$ in (\ref{eq:approximate-solution}) is composed
of approximations to nodal values $u(\vec{x}_{j})$. The FEM using
Lagrange basis functions is called \emph{Lagrange FEM} \cite[p.  36]{ciarlet2002finite}.
 In this work, FEM refers to Lagrange FEM, unless otherwise noted.

For elliptic PDEs, FEM solves (\ref{eq:mwr-eqn}) by performing integration
by parts and then substituting the boundary conditions (\ref{eq:bndy-problem-Dirichlet})
and (\ref{eq:bndy-problem-Neumann}) into the resulting integral equation.
Let $\left\langle \cdot,\cdot\right\rangle _{\Omega}$ denote the
inner product over $\Omega$,\emph{}\footnote{In functional analysis, the inner product is denoted as $(\cdot,\cdot)$.
We use $\left\langle \cdot,\cdot\right\rangle _{\Omega}$ for clarity
and for distinguishing the inner products on $\Omega$ and on boundary
$\Gamma$. See \ref{subsec:Functional-Analysis} for a review of some
relevant functional analysis concepts.} i.e., 
\begin{equation}
\left\langle \phi,\psi\right\rangle _{\Omega}=\int_{\Omega}\phi\psi\,d\vec{x},\label{eq:inner-product}
\end{equation}
which defines the $L^{2}$ norm over $\Omega$, i.e., $\Vert\phi\Vert_{L^{2}(\Omega)}=\sqrt{\langle\phi,\phi\rangle_{\Omega}}$.
Similarly, let $\left\langle \cdot,\cdot\right\rangle _{\Gamma}$
denote the inner product over $\Gamma$. For the general linear operator
$\mathcal{L}$ in (\ref{eq:2nd-order-elliptic-operator}), after integration
by parts of the first term, we obtain a variational form for each
test function $\psi_{i}$ 
\begin{align}
a(u_{h},\psi_{i}) & =\left\langle f,\psi_{i}\right\rangle _{\Omega}+\left\langle \mu\partial_{\vec{n}}u_{h},\psi_{i}\right\rangle _{\Gamma},\label{eq:weak-formulation}
\end{align}
where
\begin{equation}
a(u_{h},\psi_{i})=\int_{\Omega}\left(\vec{\nabla}\psi_{i}\cdot\left(\mu\vec{\nabla}u_{h}\right)+\vec{\nu}\cdot\vec{\nabla}u_{h}\psi_{i}+\omega^{2}u_{h}\psi_{i}\right)\,d\vec{x},\label{eq:bilinear-form-general}
\end{equation}
is the \emph{bilinear form}.   

\subsection{Prior Efforts on Unified Analysis of FDM and FEM}

The unification of the accuracy and stability analysis of FDM and
FEM has been of great interest to numerical analysts since the late
1960s \cite{fix1969fourier,strang1973analysis}. In terms of local
error analysis, this unification is straightforward by using Taylor
series, except that FDM traditionally relied on 1D Taylor series,
whereas FEM requires the higher-dimensional version. In terms of global
error analysis, Fix and Strang explored adapting Fourier analysis
from FDM to FEM on structured grids \cite{fix1969fourier}. Another
common technique used in analyzing both FDM and FEM is the Green's
functions. It was particularly successful for proving the convergence
of FDM and superconvergence of FEM in $\ell^{\infty}$ norm in 1D
or with tensor-product elements; see, e.g., \cite{douglas1974estimate,douglas1973superconvergence,leveque2007finite}.
 In this work, we unify the analysis of well-posedness of GFDM, FEM,
and AES-FEM by integrating functional analysis and approximation theory.

\section{\label{sec:Generalized-Weighted-Residuals} Generalized Weighted
Residuals}

To unify the formulations and analysis of GFDM, FEM, and AES-FEM,
we need a mathematical framework that is more general than weighted
residuals. The framework, which we refer to as \emph{generalized weighted
residuals }(\emph{GWR}), has three components: a mesh and its associated
test functions and geometric realization, a set of (generalized) Lagrange
basis functions, and a (generalized) variational form. We address
these three components using GFDM and FEM as examples, and then introduce
AES-FEM under this framework.

\subsection{Component 1: Meshes, Test Functions, and Geometry}

\subsubsection{Meshes}

Like in FEM, in GWR the domain $\Omega$ is tessellated by a \emph{mesh},
which is typically simplicial or rectangular. Without loss of generality,
we assume simplicial meshes in this work, which are triangular in
2D or tetrahedral in 3D. The triangles and tetrahedra are known as
the \emph{elements}, of which the vertices are called the \emph{nodes}.
Let $\Omega_{h}$ denote the union of the geometric realizations of
the elements, and let $\partial\Omega_{h}$ denote its boundary. Let
$\Gamma_{h}$ denote an approximation of $\Gamma$. We assume $\Gamma_{h}$
is the same as $\partial\Omega_{h}$. Let $\Omega_{h}=\Omega_{h}^{\circ}\cup\Gamma_{h,N}^{\circ}\cup\Gamma_{h,D}$,
where $\Omega_{h}^{\circ}$ denotes the interior of $\Omega_{h}$,
$\Gamma_{h,D}$ denotes the approximation to the Dirichlet boundary,
and $\Gamma_{h,N}^{\circ}=\Gamma_{h}\backslash\Gamma_{h,D}$. We refer
to the nodes in $\Omega_{h}^{\circ}$, $\Gamma_{h,N}^{\circ}$, and
$\Gamma_{h,D}$ as \emph{interior}, \emph{Neumann}, and \emph{Dirichlet}
nodes, respectively. Without loss of generality, we assume the nodes
are numbered between $1$ and $m$, where the first $n$ nodes are
those in $\Omega_{h}^{\circ}\cup\Gamma_{h,N}^{\circ}$.

Given a node $\vec{x}$, the term \textit{stencil }refers to the nodes
where the generalized Lagrange trial functions associated with $\vec{x}$
are non-zero. Note that in GFDM, stencils are referred to as \textit{stars}
\cite{jensen1972finite,perrone1975general,liszka1980finite}. See
\ref{sec:Selection-of-Stencils} for details about stencil selection
for AES-FEM.

\subsubsection{Test Functions}

In GWR, there is a \emph{test function} $\psi_{i}$ associated with
each node $\vec{x}_{i}\in\Omega_{h}$, analogous to those in (\ref{eq:mwr-eqn}).
Each test function $\psi_{i}$ has \emph{local support}, denoted by
$\Omega_{i}$, which is the closure of the subset of points in $\Omega_{h}$
such that $\psi_{i}(\vec{x})\neq0$, i.e., $\Omega_{i}=\text{cl}(\left\{ \vec{x}\mid\vec{x}\in\Omega_{h}\wedge\psi_{i}(\vec{x})\neq0\right\} )$.
Topologically, the local support is \emph{compact}, in that it contains
only a small constant number of nodes. In a Lagrange FEM, each test
function is a Lagrange basis function, such as a hat (a.k.a. pyramid)
function. Note that the test functions in FEM may be quadratic or
higher-degree polynomials, which have \emph{mid-edge}, \emph{mid-face},
and \emph{mid-cell} nodes, besides \textit{\emph{the}}\emph{ corner
}\textit{\emph{nodes}}. For FEM, the local support $\Omega_{i}$ is
composed of the union of the elements incident on $\vec{x}_{i}$.
For GFDM over an unstructured mesh, within the GWR framework, the
test function at $\vec{x}_{i}$ is the Dirac delta function at $\vec{x}_{i}$,
and the local support of a Dirac delta function is $\vec{x}_{i}$
itself. We will discuss this further in Section~\ref{subsec:Generalized-Weighted-Residuals}.
\begin{rem}
In \cite{ciarlet2002finite}, Ciarlet defined a finite element method
as a triplet: a mesh, element-based (nearly) polynomial basis, and
node-based basis functions of an $H^{1}$ space. A GWR is more general
in that the test functions may not span an $L^{2}$ or $H^{1}$ space,
which is the case in (generalized) finite difference methods.
\end{rem}

\subsubsection{\label{subsec:Geometric-Realizations}Geometric Realizations}

For numerical computations, the local support $\Omega_{i}$ must have
a \emph{geometric realization}, which is the union of the geometric
realizations of its elements. For each element $\tau$, its geometric
realization is defined through a mapping from a \emph{master element}
$e_{m}$ in the parametric space to the ``physical space'' $\mathbb{R}^{d}$.
Let $\vec{\xi}$ denote the natural coordinates in the parametric
space. Let $n_{e}$ denote the number of nodes in $e_{m}$, and let
$\vec{\xi}_{K}$ and $\vec{x}_{K}$ denote the natural coordinates
and physical coordinates, respectively, of the $K$th node in $e_{m}$,
where $1\leq K\leq n_{e}$. For example, a linear triangle has nodes
$\vec{\xi}_{1}=[0,0]$, $\vec{\xi}_{2}=[0,1]$ and $\vec{\xi}_{3}=[1,0]$.
The \emph{geometric realization} of an element $\tau$ is defined
by a Lagrange interpolation
\begin{equation}
\vec{x}_{\tau}(\vec{\xi})=\sum_{K=1}^{N_{e}}\vec{x}_{K}\varphi_{K}(\vec{\xi}),\label{eq:param-surf-element}
\end{equation}
where $K$ is the local nodal ID in $\tau$ for the $k$th node in
$\Omega_{h}$. The functions $\varphi_{K}$ are in general polynomials.
Using the interpolation theory \cite{phillips2003interpolation},
given $n_{e}$ nodes and an equal number of monomials in $\vec{\xi}$,
if the Vandermonde system in the parametric space is stable, then
the basis functions $\{\varphi_{K}\}$ are uniquely determined over
$e_{m}$. In FEM, the geometric basis functions $\{\varphi_{k}\}$
do not need to be the same as the test functions $\{\psi_{i}\}$ (and
trial functions $\{\phi_{j}\}$).

Besides the local support, we also define a ``\emph{control volume}''
$\omega_{i}$ for each node to facilitate the theoretical analysis
in Section~\ref{sec:well-posedness} by generalizing its traditional
definition. Let $\left|\omega_{i}\right|$ denote the Lebesgue measure
(namely, the area or volume) of $\omega_{i}$. The control volumes
of all the nodes partition $\Omega_{h}$, i.e., $\Omega_{h}=\cup_{i=1}^{m}\omega_{i}$
and $\left|\omega_{i}\cap\omega_{j}\right|=0$ if $i\neq j$. For
FEM, the control volume of a node can be defined by the union of its
closest points within its incident elements. For GFDM, we define the
control volume similarly or use the Voronoi cells. Note that these
control volumes are not used in computations; instead, we will use
them in analyzing well-posedness and convergence. See Figure \ref{fig:control_volume}
in \ref{sec:Selection-of-Stencils} for an example of a control volume.

\subsubsection{\label{subsec:Approximation-Power-of} Approximation Power of Lagrange
Finite Elements}

In FEM, the test (and trial) functions are Lagrange functions, which
are defined using a mapping from a parametric space $[0,1]^{d}$ to
the physical space $\mathbb{R}^{d}$, similar to that of the geometric
basis functions. More precisely, let $\{\psi_{e,J}\mid1\leq J\leq n_{e}\}$
denote the \emph{Lagrange polynomial basis} over the master element,
which satisfies the Kronecker delta property $\psi_{e,J}(\vec{\xi}_{I})=\delta_{IJ}$.
Let $\vec{x}_{\tau}(\vec{\xi}):\left[0,1\right]^{d}\supseteq e_{m}\rightarrow\tau\subset\mathbb{R}^{d}$
denote the mapping from the parametric space to the physical space
and $\vec{\xi}_{\tau}(\vec{x}):\tau\rightarrow e_{m}$ denote its
inverse mapping. A global test function $\psi_{j}(\vec{x}):\Omega\rightarrow\mathbb{R}$
is then defined as
\begin{equation}
\psi_{j}(\vec{x})=\psi_{e,J}\left(\vec{\xi}_{\tau}(\vec{x})\right)\qquad\text{if }\ensuremath{\vec{x}}\ensuremath{\in\tau},\label{eq:test-functions}
\end{equation}
where $J$ is the local ID of node $\vec{x}_{j}$ in $\tau$. By construction,
$\psi_{j}$ is a Lagrange test function over $\Omega$. A (piecewise)
smooth function $u:\Omega\rightarrow\mathbb{R}$ can be interpolated
by the basis functions over $\Omega$ by
\begin{equation}
u_{\Psi}(\vec{x})\coloneqq\sum_{j=1}^{m}u(\vec{x}_{j})\psi_{j}\left(\vec{x}\right)=\sum_{J=1}^{n_{e}}u(\vec{x}_{J})\psi_{e,J}\left(\vec{\xi}_{\tau}(\vec{x})\right),\qquad\text{if }\text{\ensuremath{\vec{x}}}\in\tau.\label{eq:interpolation-test-space}
\end{equation}

If $\varphi_{k}$ is piecewise linear, then the Lagrange test functions
$\psi_{j}(\vec{x})$ are polynomials. However, if $\varphi_{k}$ is
nonlinear, then $\psi_{j}(\vec{x})$ are no longer polynomials. Nevertheless,
$\left\Vert \vec{\nabla}^{k}u_{\Psi}-\vec{\nabla}_{\vec{x}}^{k}u\right\Vert _{\infty}$
is approximated to $\mathcal{O}(h^{p-k+1})$ within each element if
both $\left\Vert \vec{\nabla}_{\vec{\xi}}^{i}\phi_{e,J}(\vec{\xi})\right\Vert _{L^{\infty}(e_{m})}$
and $\left\Vert h^{i}\vec{\nabla}_{\vec{x}}^{i}\vec{\xi}\right\Vert _{L^{\infty}(e_{m})}$
are bounded for $i=1,\ensuremath{\dots},p+1$ \cite{ciarlet1972interpolation},
where $\vec{\nabla}^{k}$ denote the $k$th derivative tensor and
$h$ is an edge length measure.

\subsection{\label{subsec:Generalized-Lagrange-Polynomial} Component 2: Generalized
Lagrange Trial Functions}

In GWR, for each node (and more generally, at each point) in $\Omega_{h}$,
there is a set of \emph{generalized Lagrange} trial functions, which
may or may not be polynomials, and which may be interpolation (such
as in FEM) or least-squares-based (such as in GFDM).

\subsubsection{Generalized Lagrange Basis Functions}
\begin{defn}
\label{def:generalized-Lagrange} A set of functions $\left\{ \phi_{j}(\vec{x})\mid1\leq j\leq m\right\} $
form a set of \emph{generalized Lagrange basis functions} of \emph{degree-$p$}
\emph{consistency} over local support $\Omega_{i}$ with a \emph{stencil}
$\{\vec{x}_{j}\in\mathbb{R}^{d}\mid1\leq j\leq m\}$ if 
\begin{equation}
\left\Vert \sum_{j=1}^{m}u(\vec{x})\vec{\nabla}^{k}\phi_{j}(\vec{x})-\vec{\nabla}^{k}u(\vec{x})\right\Vert _{L^{\infty}(\Omega_{i})}=\left\Vert \vec{\nabla}^{p+1}u\right\Vert _{L^{\infty}(\Omega_{i})}\mathcal{O}(h^{p+1-k})\label{eq:gen-Lagrange}
\end{equation}
for a sufficiently differentiable function $u:\Omega\rightarrow\mathbb{R}$
and $k=0,\dots,p$, where $h$ is the radius of the stencil. These
basis functions are \emph{stable} over $\Omega_{i}$ if $\left\Vert h^{k}\vec{\nabla}^{k}\phi_{j}(\vec{x})\right\Vert _{L^{\infty}(\Omega_{i})}\leq C\ll\infty$
for $k=1,\dots,p$.
\end{defn}
In the above, the ``radius'' is a local length measure of the stencil,
and it can be replaced by other characteristic length measures, such
as the maximum distance between the points in the stencil. This definition
preserves two fundamental properties of degree-$p$ Lagrange basis
functions: when approximating a function $u$, the coefficient for
each $\phi_{j}$ is $u(\vec{x}_{j})$, and $u$ is approximated to
$\mathcal{O}(h^{p+1})$ consistency. In FEM, consistent Lagrange trial
(or test) functions constitute a set of generalized Lagrange basis
functions.

\subsubsection{Generalized Lagrange Polynomials}

In GFDM, the derivatives are approximated using polynomials constructed
using least squares approximations. We can express them in terms of
\emph{generalized Lagrange polynomial} (\emph{GLP}) basis functions.
\begin{defn}
\label{def:GLP-basis} Given a \emph{stencil} $\{\vec{x}_{j}\in\mathbb{R}^{d}\mid1\leq j\leq m\}$,
degree-$p$ polynomials $\left\{ \phi_{j}(\vec{x})\mid1\leq j\leq m\right\} $
form a set of \textit{generalized Lagrange polynomial }\textit{\emph{(}}\textit{GLP}\textit{\emph{)
}}\textit{basis functions} if every degree-$p$ polynomial $P(\vec{x})$
is interpolated exactly by $\sum_{j=1}^{m}P\left(\vec{x}_{j}\right)\phi_{j}(\vec{x})$.
\end{defn}
In practice, the GLP basis functions are computed from the pseudoinverse
solution of a Vandermonde system; see e.g. \cite{conley_2016}. For
completeness, we describe the procedure in \ref{sec:Computation-of-GLP}.
If the Vandermonde matrix is nonsingular, i.e., the number of points
in the stencil is equal to the number of monomials of up to degree
$p$ and the stencil is not degenerate, then a set of GLP basis functions
reduces to a set of Lagrange polynomial basis functions, which are
commonly used in the classical finite difference methods. The basis
functions are stable if the Vandermonde matrix is well conditioned.
More importantly, they are generalized Lagrange trial functions.

\begin{figure}
\begin{minipage}[t]{0.33\textwidth}%
\begin{center}
\includegraphics[width=1\textwidth]{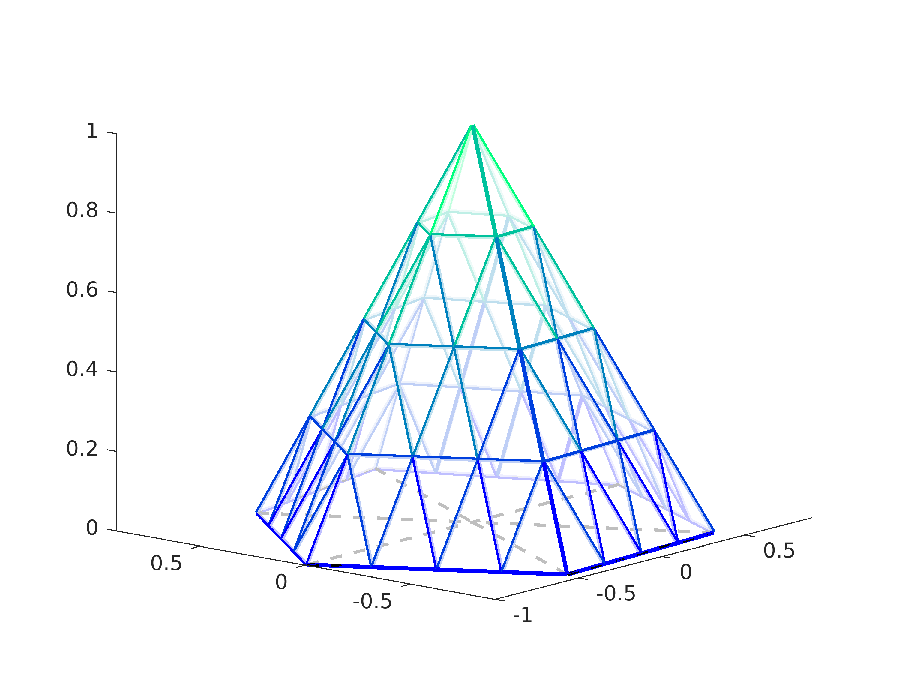}\\
{\scriptsize{}(a) A piecewise linear Lagrange}\\
{\scriptsize{}basis function in 2D.}
\par\end{center}%
\end{minipage} %
\begin{minipage}[t]{0.33\textwidth}%
\begin{center}
\includegraphics[width=1\textwidth]{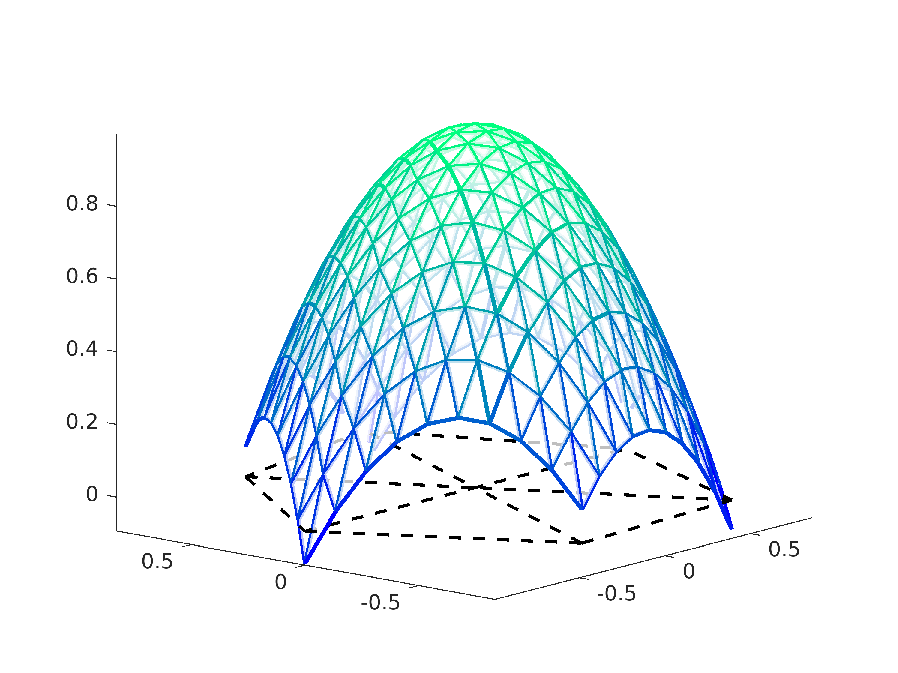}\\
{\scriptsize{}(b) A quadratic GLP}\\
{\scriptsize{}basis function in 2D.}
\par\end{center}%
\end{minipage}%
\begin{minipage}[t]{0.33\textwidth}%
\begin{center}
\includegraphics[width=1\textwidth]{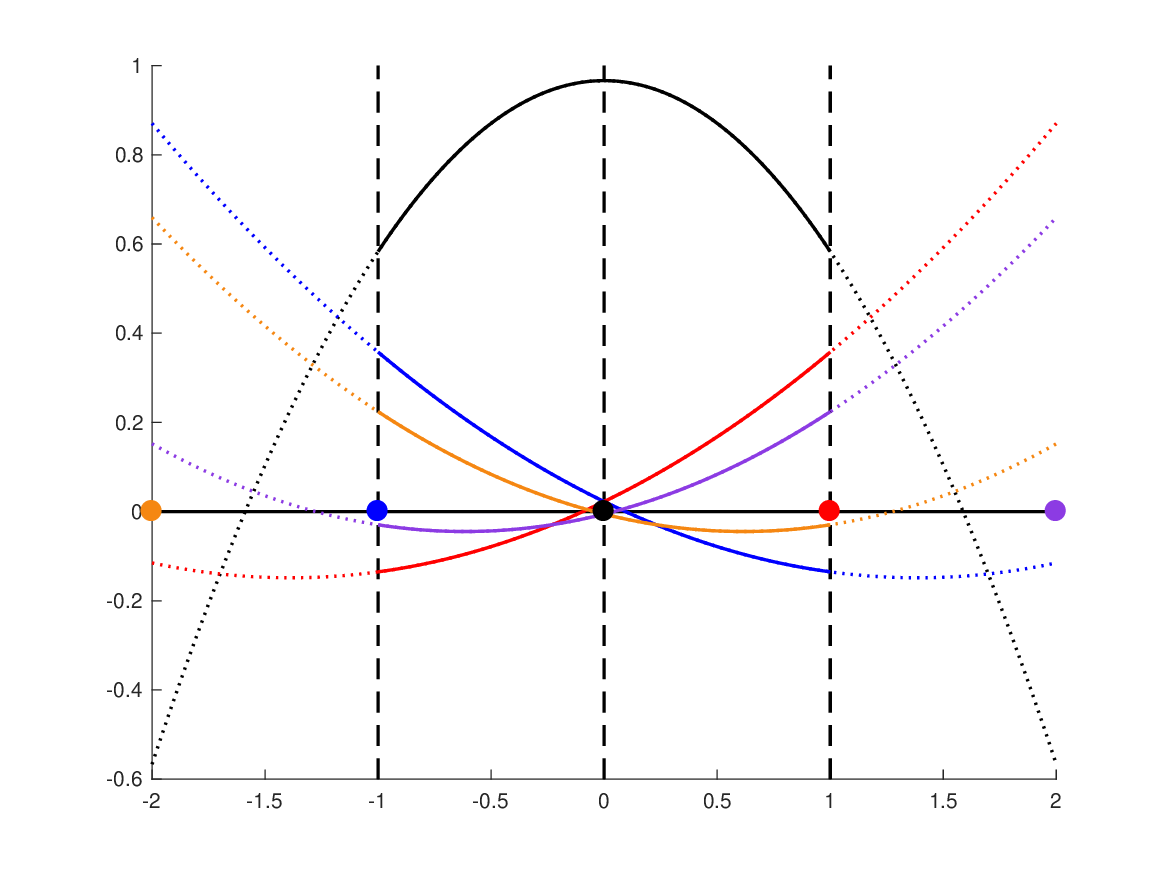}\\
{\scriptsize{}(c) A set of quadratic GLP}\\
{\scriptsize{}basis functions in 1D.}
\par\end{center}%
\end{minipage}
\raggedright{}\caption{\label{fig:Lagrange-vs-GLP}Piecewise Lagrange basis functions versus
quadratic GLP basis functions. The basis functions in (a) and (b)
are associated with the center node and have the same stencils. The
basis functions in (c) form a set of five quadratic basis functions
for the center node in 1D.}
\end{figure}

The GLP basis functions are not unique in general in that they depend
on how different points are weighted. For the AES-FEM results in Section~\ref{sec:numerical-results},
we use an inverse-distance-based weighting scheme given in (\ref{eq:weighting_scheme_AES}).
This weighting scheme tends to promote error cancellations on nearly
symmetric meshes \cite{Jiao2008} and in turn improve accuracy \cite{Li2019wls-enor}.
We defer a detailed analysis to future work.
\begin{lem}
\label{lem:deri-approximation}Given a set of stable degree-$p$ GLP
basis functions $\left\{ \phi_{j}(\vec{x})\mid1\leq j\leq m\right\} $
over $\{\vec{x}_{j}\}$, if $u$ is continuously differentiable up
to $p$th order, then (\ref{eq:gen-Lagrange}) holds.
\end{lem}
\begin{pf}
Consider the $d$-dimensional Taylor series \cite{humpherys2017foundations}
\begin{align}
u(\vec{x}_{0}+\vec{h}) & =\sum_{k=0}^{p}\frac{1}{k!}\vec{\nabla}^{k}u(\vec{x}_{0}):\vec{h}^{k}+\frac{C\left\Vert \vec{\nabla}^{p+1}u(\vec{x}_{0}+\vec{\xi})\right\Vert }{(p+1)!}\Vert\vec{h}\Vert^{p+1},\label{eq:Taylor-polynomial-residual}
\end{align}
where $\vec{h}^{k}$ denotes the $k$th tensor power of $\vec{h}$,
``:'' denotes the scalar product of $k$th-order tensors, $\Vert\vec{\xi}\Vert\leq\Vert\vec{h}\Vert$,
and $\vert C\vert\leq1$. Let $u_{p}$ denote the degree-$p$ Taylor
polynomial (i.e., the first term in (\ref{eq:Taylor-polynomial-residual})).
By definition, $u_{p}=\sum_{j=1}^{m}u_{p}\left(\vec{x}_{j}\right)\phi_{j}(\vec{x})$.
Let $\delta u=u_{p}-u$. 
\[
\left\Vert \sum_{j=1}^{m}u\left(\vec{x}_{j}\right)\vec{\nabla}^{k}\phi_{j}-\vec{\nabla}^{k}u\right\Vert _{\infty}\leq\left\Vert \sum_{j=1}^{m}\delta u\left(\vec{x}_{j}\right)\vec{\nabla}^{k}\phi_{j}\right\Vert _{\infty}+\left\Vert \vec{\nabla}^{k}\delta u\right\Vert _{\infty},
\]
where both terms are bounded by $\left\Vert \vec{\nabla}^{p+1}u(\vec{x})\right\Vert _{\infty}\mathcal{O}\left(h^{p-k+1}\right)$.
\end{pf}
We note that there are some differences between the Lagrange basis
functions in FEM and the GLP basis functions. First, the GLP basis
functions are least-squares based, so they, in general, do not satisfy
the Kronecker delta property, i.e., $\phi_{j}(\vec{x}_{i})\neq\delta_{ij}$.
Second, the Lagrange basis functions in FEM are \emph{$C^{0}$ continuous},
whereas the GLP basis functions are \emph{quasicontinuous} in that
they are smooth over the local support, but they do not vanish exactly
along the boundary. We illustrate these differences in Figure~\ref{fig:Lagrange-vs-GLP},
which shows (a) a 2D FEM hat function, (b) a 2D quadratic GLP basis
function at a node over the same stencil, and (c) a set of GLP basis
functions at a node in 1D. Third, Lagrange basis functions in FEM
are defined based on a mapping between the reference domain to the
physical domain, and hence they depend on element shapes and are in
general not polynomials with nonlinear geometric realizations, whereas
the GLP basis functions do not depend on the element shapes and are
true polynomials. Finally, the GLP basis functions are defined locally
at a node (or a center point), and they do not necessarily define
a global set of trial functions over $\Omega_{h}$. For these reasons,
GWR requires a more general variational form than that used in FEM.

\subsection{\label{subsec:Generalized-Weighted-Residuals} Component 3: Generalized
Variational Form}

Consider a node $\vec{x}_{i}\in\Omega_{h}^{\circ}\cup\Gamma_{h,N}^{\circ}$,
and let $\psi_{i}$ denote the test function associated with $\vec{x}_{i}$.
Let $\vec{\Phi}_{i}$ denote the vector containing the trial functions
$\{\phi_{ij}\}$ associated with the node $\vec{x}_{i}$. Let $\vec{\Phi}_{i,\circ}$
and $\vec{\Phi}_{i,D}$ denote the subvectors of $\vec{\Phi}_{i}$
corresponding to the nodes in $\Omega_{h}^{\circ}\cup\Gamma_{h,N}^{\circ}$
and $\Gamma_{h,D}$, respectively, i.e., $\vec{\Phi}_{i,\circ}=\vec{\Phi}_{i,1:n}$
and $\vec{\Phi}_{i,D}=\vec{\Phi}_{i,n+1:m}$ in MATLAB-style colon-notation.
Let $\vec{u}_{h}$ and $\vec{u}_{D}$ be composed of the nodal values
associated with $\vec{\Phi}_{i,\circ}$ and $\vec{\Phi}_{i,D}$, respectively.
In GWR, we approximate the solution of (\ref{eq:bndy-value-problem})
locally about $\vec{x}_{i}$ by 
\[
u_{h,i}=\vec{\Phi}_{i}^{T}\begin{bmatrix}\vec{u}_{h}\\
\vec{u}_{D}
\end{bmatrix}=\vec{\Phi}_{i,\circ}^{T}\vec{u}_{h}+\vec{\Phi}_{i,D}^{T}\vec{u}_{D}.
\]
Let us first consider a BVP with Dirichlet boundary conditions. We
define a \emph{generalized variational form} (GVF) corresponding to
(\ref{eq:bndy-value-problem}) for an interior node $\vec{x}_{i}\in\Omega_{h}^{\circ}$
as 
\begin{align}
a_{\circ}\left(u_{h,i},\psi_{i}\right) & =\langle f,\psi_{i}\rangle_{\Omega_{h}},\label{eq:varform-gwr-interior}
\end{align}
where 
\begin{equation}
a_{\circ}\left(u_{h,i},\psi_{i}\right)=\left\langle \mathcal{L}u_{h,i},\psi_{i}\right\rangle _{\Omega_{h}},\label{eq:bilinear-form-gwr-interior}
\end{equation}
is the \emph{generalized bilinear form} associated with $\psi_{i}$.
Here, the inner product is computed over the local support $\Omega_{i}$
of $\psi_{i}$. These forms are ``generalized'' in that $\psi_{i}$
may be a generalized function (such as a Dirac delta function). Note
that if $\psi_{i}$ vanishes along $\Gamma$, the bilinear form in
(\ref{eq:bilinear-form-general}) and the generalized bilinear form
(\ref{eq:bilinear-form-gwr-interior}) are mathematically equivalent
to each other due to Green's identities. Computationally, however,
(\ref{eq:bilinear-form-general}) requires $\psi_{i}$ to be at least
$C^{0}$ and (\ref{eq:bilinear-form-gwr-interior}) requires $u_{h,i}$
to be at least $C^{1}$.

\subsubsection{Generalized Variational Form of GFDM}

In (G)FDM, since the test functions are Dirac delta functions, which
are not $C^{0}$, we must use (\ref{eq:bilinear-form-gwr-interior})
instead of (\ref{eq:bilinear-form-general}) computationally. Neumann
boundary conditions can be incorporated into the generalized variational
form; see \cite[p. 31]{leveque2007finite} for derivations of second-order
FDM in 1D. Alternatively, one-sided differences can be used to convert
(\ref{eq:bndy-problem-Neumann}) directly into an algebraic equation
at each Neumann node \cite[p. 32]{leveque2007finite}. 

\subsubsection{Variational Form of Lagrange FEM}

For a Lagrange FEM with strongly imposed Dirichlet boundary conditions,
its variational form is also a GVF. Since the trial functions are
$C^{0}$, we must use the bilinear form (\ref{eq:bilinear-form-general})
computationally for $\vec{x}_{i}\in\Omega_{h}^{\circ}$. For BVP with
Neumann boundary conditions, assuming accurate geometry, one can substitute
$g_{N}\approx\mu\partial_{\vec{n}}u_{h}$ into (\ref{eq:weak-formulation})
for each node $\vec{x}_{i}$ in $\Gamma_{h,N}$, which \emph{weakly}
imposes Neumann boundary conditions (\ref{eq:bndy-problem-Neumann}).

\subsection{\label{subsec:AES-FEM}Adaptive Extended Stencil FEM}

\emph{Adaptive Extended Stencil FEM }or\emph{ AES-FEM}, is a GWR method
that combines the features of GFDM and FEM. In particular, AES-FEM
uses a piecewise linear FEM mesh, and its associated hat functions
as test functions in the interior. Its trial functions are the GLP
trial functions. We say that an AES-FEM is degree-$p$ if its trial
functions are degree-$p$ GLP basis functions at each node (or more
generally, at each point). Because of the least-squares nature of
GLP basis functions, the number of nodes in the stencil can be chosen
adaptively to ensure the well-conditioning of the Vandermonde matrix.
This is the reason for the name ``adaptive extended stencil'' in
AES-FEM; we describe the selection of the stencils and its adaptation
in \ref{sec:Selection-of-Stencils}. Because the test functions are
$C^{0}$ and the trial functions are differentiable to $p$th order,
we can use either (\ref{eq:bilinear-form-general}) or (\ref{eq:bilinear-form-gwr-interior})
computationally, which would give the same results up to machine precision.

In \cite{conley_2016}, the authors considered quadratic AES-FEM and
showed its second-order accuracy with Dirichlet boundary conditions.
The focus of this work is to establish the well-posedness and convergence
of higher-degree AES-FEM for Neumann boundary conditions over polygonal
domains. As in FEM, we can substitute $g_{N}=\mu\partial_{\vec{n}}u_{h}$
into (\ref{eq:weak-formulation}) for each node $\vec{x}_{i}$ in
$\Gamma_{h,N}$ for AES-FEM. For curved domains, Neumann boundary
conditions can be imposed using higher-order boundary representations
similar to \cite{sevilla2008nurbs}, as described in Section~\ref{subsec:Approximation-Power-of},
or using Dirac delta test functions and one-sided differences as in
GFDM. In this work, we assume polygonal domains and impose Neumann
boundary conditions weakly as in FEM. We defer the analysis and comparison
of different techniques for imposing Neumann conditions over curved
boundaries to future work.
\begin{rem}
GFDM requires a cloud of nodes, rather than a mesh with elements,
so it is considered a meshless method. However, to construct the stencils
(a.k.a. stars), one needs a data structure, such as a quadtree/octree
\cite{gavete2015approach}. Because AES-FEM involves integration (see
\ref{sec:Overview-of-AES-FEM}), we use a mesh and an associated data
structure to construct the stencils and to provide the elements for
integration. The use of integration allows AES-FEM to enjoy additional
error cancellations and hence better accuracy than GFDM on nearly
symmetric meshes, as demonstrated in Section \ref{sec:numerical-results}.
\end{rem}

\section{\label{sec:well-posedness} Well-Posedness in $\ell^{p}$ Norm}

We first establish the \emph{well-posedness}, and more precisely,
\emph{algebraic invertibility},\emph{ }of generalized weighted residual
methods for elliptic BVPs. Like that of FEM \cite{ern2013theory},
this invertibility implies the existence and uniqueness of a solution
in exact arithmetic. We focus on AES-FEM while making the analysis
general enough for GFDM. In Section~\ref{sec:convergence}, we will
establish the convergence of AES-FEM in the presence of rounding errors.

\subsection{\label{subsec:Algebraic-Equations-of} Algebraic Equations of GWR}

To analyze a GWR method, we first convert its GVF into a system of
linear equations. For generality, we assume the bilinear form (\ref{eq:bilinear-form-general}),
with strongly imposed Dirichlet boundary conditions and weakly imposed
Neumann boundary conditions. For GFDM, we assume only Dirichlet boundary
conditions, and the GVF is evaluated with (\ref{eq:bilinear-form-gwr-interior})
instead of (\ref{eq:bilinear-form-general}).

From (\ref{eq:bilinear-form-general}), we obtain an $n\times n$
linear system in $\vec{u}_{h}$, namely, 
\begin{equation}
\vec{A}\vec{u}_{h}=\vec{b},\label{eq:linear-system}
\end{equation}
where the $i$th row of $\vec{A}$ and $\vec{b}$ are
\begin{align}
\vec{a}_{i}^{T} & =a(\vec{\Phi}_{i,\circ}^{T},\psi_{i}),\label{eq:coefficient-matrix-gwr}\\
b_{i} & =\langle f,\psi_{i}\rangle_{\Omega_{h}}-a\left(\vec{\Phi}_{i,D}^{T}\vec{u}_{D},\psi_{i}\right)+\left\langle \mu\partial_{\vec{n}}u_{h},\psi_{i}\right\rangle _{\Gamma_{h,N}},\label{eq:load-vector-gwr}
\end{align}
respectively. For the Poisson equation, $\vec{A}$ is the \emph{stiffness
matrix} and $\vec{b}$ is the \emph{load vector} in FEM. For generality,
we simply refer to $\vec{A}$ as the \emph{coefficient matrix}. Another
important matrix is the \emph{mass matrix}, which is the coefficient
matrix of the \emph{constrained projection} associated with the GWR,
\begin{equation}
\vec{M}\vec{u}_{P}=\vec{b}_{P},\label{eq:constrained-projection}
\end{equation}
where the $i$th row of $\vec{M}$ and $\vec{b}$ are
\begin{align*}
\vec{m}_{i}^{T} & =\langle\vec{\Phi}_{i,\circ}^{T},\psi_{i}\rangle_{\Omega_{h}},\\
b_{P,i} & =\left\langle u-u_{D}\vec{\Phi}_{i,D}^{T}\vec{u}_{D},\psi_{i}\right\rangle _{\Omega_{h}}.
\end{align*}
For Galerkin FEM, the constrained projection is known as the \emph{$L^{2}$
projection} \cite[p. 132]{ern2013theory}. Note that Neumann constraints
are not imposed explicitly in the constrained projection because they
are satisfied weakly automatically.

\subsection{\label{subsec:Algebraic-Error-Analysis} Algebraic Error Analysis
in $\ell^{p}$ Norm}

To analyze the solutions of (\ref{eq:linear-system}), we consider
the nodal values in $\ell^{p}$ norm. Given a vector $\vec{v}\in\mathbb{R}^{n}$,
its $\ell^{p}$ norm is $\Vert\vec{v}\Vert_{p}=\sqrt[^{p}]{\sum_{i=1}^{n}\left|v_{i}^{p}\right|}$
, where $1\leq p\leq\infty$. Note that this $p$ is independent of,
and different from, the degree of the basis functions. We will primarily
use the $\ell^{2}$ or $\ell^{\infty}$ norm (i.e., $p=2$ or $p=\infty$).
The $\ell^{p}$\emph{ norm} of $\vec{A}\in\mathbb{R}^{n\times n}$
is $\left\Vert \vec{A}\right\Vert _{p}=\sup_{\left\Vert \vec{v}\right\Vert _{p}=1}\left\Vert \vec{A}\vec{v}\right\Vert _{p}$.

Let $\vec{u}_{h}\in\mathbb{R}^{n}$ denote the solution vector of
a GWR method on a mesh with $n$ nodes. Let $\vec{u}_{I}$ denote
the vector composed of $u(\vec{x}_{i})$. The \emph{error vector}
is
\begin{equation}
\delta\vec{u}=\vec{u}_{h}-\vec{u}_{I}.\label{eq:error-vector}
\end{equation}
Consider the linear system (\ref{eq:linear-system}). Its \emph{residual
vector} is 
\begin{equation}
\vec{r}=\vec{b}-\vec{A}\vec{u}_{I}=\vec{A}(\vec{u}_{h}-\vec{u}_{I})=\vec{A}\delta\vec{u}.\label{eq:residual}
\end{equation}
Suppose $\vec{A}\in\mathbb{R}^{n\times n}$ is nonsingular, and assume
exact arithmetic. Then,
\begin{equation}
\left\Vert \delta\vec{u}\right\Vert _{p}\leq\left\Vert \vec{A}^{-1}\right\Vert _{p}\left\Vert \vec{r}\right\Vert _{p},\label{eq:backward-error}
\end{equation}
where $\left\Vert \vec{A}^{-1}\right\Vert _{p}$ is the \emph{absolute
condition number} \emph{in $p$-norm} of the linear system (\ref{eq:linear-system}).
Note that given $\vec{A}\in\mathbb{R}^{n\times n}$ and $\vec{B}\in\mathbb{R}^{n\times n}$,
$\left\Vert \vec{A}\vec{B}\right\Vert _{p}\leq\left\Vert \vec{A}\right\Vert _{p}\left\Vert \vec{B}\right\Vert _{p}$.

\subsection{\label{subsec:Invertibility-of-Variational} Well-Posedness of Constrained
Projection in $\ell^{p}$ Norm}

We first apply backward error analysis to the constrained projection
(\ref{eq:constrained-projection}). It is an important base case for
the analysis of elliptic PDEs. In particular, consider a perturbation
$\delta u$ to $u$ in the right-hand side of (\ref{eq:constrained-projection}).
This leads to a perturbation $\delta\vec{b}$ in $\vec{b}_{P}$. From
(\ref{eq:backward-error}), the perturbation $\delta\vec{u}$ in $\vec{u}_{P}$
is bounded by
\begin{equation}
\left\Vert \delta\vec{u}\right\Vert _{p}\leq\left\Vert \vec{M}^{-1}\right\Vert _{p}\left\Vert \delta\vec{b}\right\Vert _{p}.\label{eq:perturbation-inequality}
\end{equation}
If $\delta u$ is $C^{0}$ continuous, $\left\Vert \delta\vec{b}\right\Vert _{p}=(1+\mathcal{O}(h))\left\Vert \vec{M}\left\llbracket \delta u(\vec{x}_{i})\right\rrbracket \right\Vert _{p}\leq\mathcal{O}(1)\left\Vert \vec{M}\right\Vert _{p}\left\Vert \left\llbracket \delta u(\vec{x}_{i})\right\rrbracket \right\Vert _{p}$,
so $\left\Vert \delta\vec{u}\right\Vert _{p}=\mathcal{O}(1)\left\Vert \left\llbracket \delta u(\vec{x}_{i})\right\rrbracket \right\Vert _{p}$
if $\kappa_{p}(\vec{M})=\mathcal{O}(1)$; on the other hand, if $\kappa_{p}(\vec{M})=\mathcal{O}(h^{-\alpha})$
for some $\alpha>0$, an $\mathcal{O}(1)$ continuous perturbation
$\delta u$ in the right-hand side of (\ref{eq:constrained-projection})
may lead to an $\mathcal{O}(h^{-\alpha})$ perturbation in $\delta\vec{u}$
in $\ell^{p}$ norm. Hence, to be consistent with the classical Hadamard's
notion of well-posedness of variational methods \cite[p. 82]{ern2013theory},
we define a well-posed constrained projection as follows.
\begin{defn}
\label{def:well-conditioned-constrained-proj} A constrained projection
is \emph{well-posed} in $\ell^{p}$ norm for $1\leq p\leq\infty$
independently of $h$ if it is \emph{well-conditioned}, i.e.,
\begin{equation}
\kappa_{p}(\vec{M})=\left\Vert \vec{M}\right\Vert _{p}\left\Vert \vec{M}^{-1}\right\Vert _{p}=\mathcal{O}(1).\label{eq:quasiuniformity-1}
\end{equation}
\end{defn}
In practice, the well-posedness requires \emph{quasiuniform meshes}.
\begin{defn}
\label{def:quasiuniform-mesh} A type of GWR meshes is \emph{quasiuniform
}if the ratio of the largest and smallest control volumes of the test
functions is bounded independently of mesh resolution, i.e., $\sup_{i}\left|\omega_{i}\right|/\inf_{i}\left|\omega_{i}\right|=\mathcal{O}(1)$.
\end{defn}
For FEM, Definition~\ref{def:quasiuniform-mesh} is satisfied with
the classical definition of quasiuniform meshes (e.g., \cite{brenner2008mathematical,strang1973analysis}),
which requires the ratio between the largest and smallest elements
to be bounded. Definition~\ref{def:quasiuniform-mesh} is more general
and also applies to GFDM, of which the control volumes have nonzero
measures, but the local support of a Dirac delta function has a zero
measure; see Remark~\ref{rem:GFDM-cons-proj}. For AES-FEM, we note
the following fact.
\begin{thm}
\label{thm:well-posedness} Constrained projection by AES-FEM is well-posed
in $\ell^{p}$ norm for $1\leq p\leq\infty$ on a sufficiently fine
quasiuniform mesh with consistent and stable GLP basis functions and
test functions.
\end{thm}
This theorem is similar to that of the well-posedness of $L^{2}$
projections \cite[p. 387]{ern2013theory}, but there are two complications.
First, the GLP basis functions are not global basis functions over
$\Omega_{h}$, so we cannot use functional analysis directly. Second,
Theorem~\ref{thm:well-posedness} is not limited to $\ell^{2}$ norms,
so we cannot use eigenvalue analysis either. To address the first
issue, we define a global basis function by blending the GLP basis
functions using $\{\psi_{i}\}$ to obtain a $C^{0}$ basis function,
i.e., $\hat{\vec{\Phi}}=\sum_{i=1}^{m}\psi_{i}\vec{\Phi}_{i}$, which
is composed of
\begin{equation}
\hat{\phi}_{j}=\sum_{i=1}^{n}\psi_{i}\phi_{ij}.\label{eq:case-fem-basis-functions}
\end{equation}
These blended basis functions have the same approximation order as
the GLP basis functions \cite{jiao2012reconstructing}, and hence
\[
\left|\hat{\vec{v}}^{T}\vec{M}\hat{\vec{u}}\right|=\left\langle \hat{\vec{\Phi}}^{T}\hat{\vec{u}},\vec{\Psi}^{T}\hat{\vec{v}}\right\rangle \left(1+\mathcal{O}\left(h^{p+1}\right)\right).
\]
To overcome the second issue, we make use of Singer's representation
theorem \cite{singer1957linear}. We omit the detailed proof, which
is similar to that of Theorem~\ref{thm:well-posed-elliptic-bvp}
below.

\subsection{\label{subsec:Conditions-of-Algebraic} Well-Posedness for Elliptic
BVPs}

Consider a perturbation $\delta f$ to $f$ in the right-hand side
of (\ref{eq:constrained-projection}). This leads to a perturbation
$\delta\vec{b}$ in $\vec{b}_{P}$, where $n^{-\frac{1}{p}}\left\Vert \delta\vec{b}\right\Vert _{p}\leq\left\Vert \vec{M}\right\Vert _{p}\left\Vert \delta f\right\Vert _{\infty}$.
From (\ref{eq:backward-error}), the perturbation $\delta\vec{u}$
in $\vec{u}_{h}$ is bounded by
\[
\left\Vert \delta\vec{u}\right\Vert _{p}\leq\left\Vert \vec{A}^{-1}\right\Vert _{p}\left\Vert \delta\vec{b}\right\Vert _{p}.
\]
Hence, $n^{-\frac{1}{p}}\left\Vert \delta\vec{u}\right\Vert _{p}=\mathcal{O}(1)\left\Vert \delta f\right\Vert _{\infty}$
if $\left\Vert \vec{A}^{-1}\right\Vert _{p}\left\Vert \vec{M}\right\Vert _{p}=\mathcal{O}(1)$.
On the other hand, if $\left\Vert \vec{A}^{-1}\right\Vert _{p}\left\Vert \vec{M}\right\Vert _{p}=\mathcal{O}(h^{-\alpha})$
for some $\alpha>0$, an $\mathcal{O}(1)$ $C^{0}$ continuous perturbation
$\delta f$ in the right-hand side of (\ref{eq:constrained-projection})
may lead to an $\mathcal{O}(h^{-\alpha})$ perturbation in $\delta\vec{u}$
in $\vec{u}_{h}$. Hence, we define well-posedness as follows.
\begin{defn}
\label{def:alg-inv} Given an elliptic BVP, let $\vec{A}$ and $\vec{M}$
denote the coefficient and mass matrices defined in (\ref{eq:linear-system})
and (\ref{eq:constrained-projection}), respectively. A GWR method
is \emph{well-posed in $\ell^{p}$ norm} for $1\leq p\leq\infty$
independently of $h$ if
\begin{equation}
\left\Vert \vec{A}^{-1}\vec{M}\right\Vert _{p}\leq\left\Vert \vec{A}^{-1}\right\Vert _{p}\left\Vert \vec{M}\right\Vert _{p}=\mathcal{O}(1).\label{eq:invertibility}
\end{equation}
\end{defn}
For AES-FEM, we note the following theorem.
\begin{thm}
\label{thm:well-posed-elliptic-bvp} Given an elliptic BVP, AES-FEM
is well-posed in $\ell^{p}$ norm for $1\leq p\leq\infty$ on a sufficiently
fine quasiuniform mesh with consistent and stable GLP basis functions
and test functions.
\end{thm}
Similar to Theorem~\ref{thm:well-posedness}, the proof requires
an adaptation by using $\hat{\vec{\Phi}}$ to construct a $C^{0}$
approximation in order to apply Singer's representation theorem. Similar
to the Lax-Milgram lemma, the proof also involves an assumption of
invertibility of the PDE in infinite dimensions, and a boundedness
assumption due to Friedrichs' inequality \cite[p. 104]{brenner2008mathematical},
which is more general than the Poincaré inequality \cite[p. 489]{ern2013theory}.
For completeness, we give the proof as follows.
\begin{pf}
Let $\hat{\vec{u}}=\arg\inf_{\left\Vert \vec{u}\right\Vert _{p}=1}\left\Vert \vec{A}\vec{u}\right\Vert _{p}$
and $\hat{u}=\hat{\vec{\Phi}}^{T}\hat{\vec{u}}$. Let $p'$ denote
\emph{Hölder's conjugate} of $p$, i.e., $1/p+1/p'=1$, and let $\hat{v}=\arg\sup_{\left\Vert v\right\Vert _{L^{p'}}=1}a\left(\hat{u},v\right)_{\Omega}$.
Since $\hat{u}$ is $C^{0}$ continuous, due to Singer's representation
theorem, there exists a solution $\hat{\vec{v}}$ such that $a\left(\hat{u},\vec{\Psi}^{T}\hat{\vec{v}}\right)=a\left(\hat{u},\hat{v}\right)$.
If the PDE is invertible in infinite dimensions,
\begin{equation}
\exists C,\qquad\inf_{\left\Vert v\right\Vert _{L^{p}(\Omega)}=1}\sup_{\left\Vert v\right\Vert _{L^{p'}(\Omega)}=1}a(u,v)\geq C>0.\label{eq:inf-sup-cond}
\end{equation}
Due to Friedrichs' inequality, $\left\Vert \hat{u}\right\Vert _{L^{p}}\left\Vert \vec{\Psi}^{T}\hat{\vec{v}}\right\Vert _{L^{p'}}\leq\Theta(1)a\left(\hat{u},\hat{v}\right)$
and $\left\Vert \vec{\Psi}^{T}\hat{\vec{v}}\right\Vert _{L^{p'}}=\Theta(1)$.
Furthermore, $\left\Vert \vec{A}^{-1}\right\Vert _{p}^{-1}=\sup_{\vec{v}\neq\vec{0}}\left|\vec{v}^{T}\vec{A}\hat{\vec{u}}\right|/\left\Vert \vec{v}\right\Vert _{p'}$,
Hence,
\[
\left\Vert \hat{\vec{v}}\right\Vert _{p'}\left\Vert \vec{A}^{-1}\right\Vert _{p}^{-1}\geq\left|\hat{\vec{v}}^{T}\vec{A}\hat{\vec{u}}\right|=a\left(\hat{u},\hat{v}\right)(1+\mathcal{O}(h))\geq(C+\mathcal{O}(h))\left\Vert \hat{u}\right\Vert _{L^{p}}.
\]
On a quasiuniform mesh, 
\[
\left\Vert \hat{u}\right\Vert _{L^{p}}=\left(\sum_{i=1}^{m}\int_{\omega_{i}}\left|\hat{u}\right|^{p}\,d\vec{x}\right)^{1/p}=\Theta(1)\left(\max\left|\omega_{i}\right|\right)^{1/p}.
\]
Since $\left\Vert \hat{v}\right\Vert _{L^{p'}}=1$, 
\[
\left\Vert \hat{\vec{v}}\right\Vert _{p'}=\left(\sum_{i=1}^{m}\left|\hat{v}_{i}\right|^{p'}\right)^{1/p'}\leq\frac{\left(\sum_{i=1}^{m}\left|\hat{v}_{i}^{p'}\left|\omega_{i}\right|\right|\right)^{1/p'}}{\left(\min\left|\omega_{i}\right|\right)^{1/p'}}=\frac{\left\Vert \vec{\Psi}^{T}\hat{\vec{v}}\right\Vert _{L^{p'}}+\mathcal{O}(h)}{\left(\min\left|\omega_{i}\right|\right)^{1/p'}}.
\]
Since $\left\Vert \vec{M}\right\Vert _{p}=\max\left|\omega_{i}\right|\Theta(1)=\min\left|\omega_{i}\right|\Theta(1)$,
\[
\left\Vert \hat{\vec{v}}\right\Vert _{p'}/\left\Vert \hat{u}\right\Vert _{L^{p}}\leq\Theta(1)\left(\min\left|\omega_{i}\right|\right)^{-1/p-1/p'}=\Theta(1)\left\Vert \vec{M}\right\Vert _{p}^{-1},
\]
and 
\[
\left\Vert \vec{A}^{-1}\vec{M}\right\Vert _{p}\leq\left\Vert \vec{M}\right\Vert _{p}\left\Vert \vec{A}^{-1}\right\Vert _{p}=\frac{\left\Vert \vec{M}\right\Vert _{p}\left\Vert \hat{\vec{v}}\right\Vert _{p'}}{\left|\hat{\vec{u}}\vec{A}\hat{\vec{v}}\right|}\leq\frac{\left\Vert \vec{M}\right\Vert _{p}\left\Vert \hat{\vec{v}}\right\Vert _{p'}}{(C+\mathcal{O}(h))\left\Vert \hat{u}\right\Vert _{L^{p}}}=\mathcal{O}(1).
\]
\end{pf}
Let $\delta\vec{u}_{h-I}=\vec{u}_{h}-\vec{u}_{I}$, where $\vec{u}_{h}$
is composed of nodal solutions of AES-FEM and $\vec{u}_{I}$ is composed
of the interpolated nodal values. It is easy to see that the above
proof also applies to FEM, simply by replacing $\hat{\vec{\Phi}}$
with the Lagrange basis functions $\vec{\Phi}$ in the proof. Hence,
Theorems~\ref{thm:well-posedness} and \ref{thm:well-posed-elliptic-bvp}
both apply to FEM.
\begin{rem}
\label{rem:GFDM-cons-proj} We can generalize Theorems~\ref{thm:well-posedness}
and \ref{thm:well-posed-elliptic-bvp} to GFDM as follows. Let $\hat{\vec{M}}=\vec{W}\vec{M}$
and $\hat{\vec{A}}=\vec{W}\vec{A}$, where $\vec{W}$ is a diagonal
matrix with $w_{i}=\left|\omega_{i}\right|$. Let $\hat{\vec{\Psi}}$
denote the vector of hat functions over the mesh. Then, it is easy
to show that $\left|\hat{\vec{v}}^{T}\hat{\vec{M}}\hat{\vec{u}}\right|=\left\langle \hat{\vec{\Phi}}^{T}\hat{\vec{u}},\hat{\vec{\Psi}}^{T}\hat{\vec{v}}\right\rangle \left(1+\mathcal{O}\left(h\right)\right)$
and $\left|\hat{\vec{v}}^{T}\hat{\vec{A}}\hat{\vec{u}}\right|=a\left(\hat{\vec{\Phi}}^{T}\hat{\vec{u}},\hat{\vec{\Psi}}^{T}\hat{\vec{v}}\right)\left(1+\mathcal{O}\left(h\right)\right)$.
By replacing $\vec{A}$ with $\hat{\vec{A}}$, the proof for Theorem~\ref{thm:well-posed-elliptic-bvp}
applies to GFDM with Dirichlet boundary conditions; similarly for
Theorem~\ref{thm:well-posedness}. 
\end{rem}
Note that $\Vert\vec{A}\Vert_{p}=\Vert\vec{M}\Vert_{p}\mathcal{O}(h^{-2})$
on a quasiuniform mesh with stable GLP basis functions. Hence a corollary
of Theorem~\ref{thm:well-posed-elliptic-bvp} is that
\begin{equation}
\kappa_{p}(\vec{A})=\left\Vert \vec{A}\right\Vert _{p}\left\Vert \vec{A}^{-1}\right\Vert _{p}\leq\left\Vert \vec{M}\right\Vert _{p}\left\Vert \vec{A}^{-1}\right\Vert _{p}\mathcal{O}(h^{-2})=\mathcal{O}(h^{-2})\label{eq:matrix-condition-number}
\end{equation}
on a quasiuniform mesh. Similarly, $\kappa_{p}(\vec{M}^{-1}\vec{A})=\mathcal{O}(h^{-2})$.
This condition number estimation is well known in $\ell^{2}$ norm
for FDM and FEM (see e.g. \cite{ern2013theory,leveque2007finite}).

\subsection{\label{subsec:Mesh-Dependency} Mesh Dependency for Well-Posedness}

From the preceding analysis, it is clear that all the GWR methods
have some level of dependency on meshes. In particular, all the methods
require the \emph{quasiuniformity} of control volumes of the nodes.
For FEM and AES-FEM, this is equivalent to the quasiuniformity of
the local support of the test functions. For GFDM, although the computation
does not depend on a mesh, the quasiuniformity imposes restrictions
on the distributions of the nodes.

Besides quasiuniformity, FEM requires \emph{well-shaped} elements,
because the Lagrange trial and test functions are based on the transformation
from the parametric space to the physical space. For linear elements
over polytopal domains, the well-shapedness requires the angles within
the elements to be bounded away from $\pi$ and 0 \cite{babuska1976angle,Shewchuk02whatis},
which is needed for the stability of interpolations and derivative
approximations. If some elements contain angles that are too small,
the stiffness matrix in FEM may become ill-conditioned \cite{shewchuk2002delaunay}.
For high-order elements, the nodes must be well-positioned within
the master elements so that the Lagrange basis functions are stable.
In contrast, the well-posedness of GFDM does not depend on element
shapes, but the stability of the GLP basis functions does depend on
the selection of the stencils. Similarly, AES-FEM also depends on
the stencils for the stability of its trial functions. However, the
stencils in GFDM and AES-FEM can be adapted more easily due to their
least squares nature. If the generalized bilinear form (\ref{eq:bilinear-form-gwr-interior})
is used, then AES-FEM also depends on the stability of the Lagrange
test functions in the parametric space, but it does not require well-shaped
elements in the physical space. If the bilinear form (\ref{eq:bilinear-form-general})
is used over linear elements, then (\ref{eq:bilinear-form-general})
is equal to (\ref{eq:bilinear-form-gwr-interior}) to machine precision,
and hence there is no dependency on element shape either. In addition,
high-order AES-FEM requires only first-order meshes for its implementation
at least in the interior of the domain, so its implementation is simpler
than that of high-order finite elements, which requires high-order meshes.

For geometries with curved boundaries, Lagrange FEM typically uses
isoparametric elements, for which the well-posedness depends on the
Ciarlet-Raviart condition \cite{ciarlet1972interpolation}. Assuming
stable Lagrange basis functions, the Ciarlet-Raviart condition requires
the $k$th derivatives for the mapping from the physical space to
the parametric space to be bounded for $k=2,\dots,p+1$, i.e., $h^{k}\left\Vert \vec{\nabla}_{\vec{x}}^{k}\vec{\xi}\right\Vert _{\infty}\leq C\ll\infty$
for $\vec{\xi}$ in the master element. Mathematically, this condition
is needed due to the high-order chain rule, also known as the Faà
di Bruno’s formula \cite{ma2009higher},
\[
\partial_{\vec{d}}^{k}\phi\left(\vec{\xi}(\vec{x}_{0})\right)=\sum_{m_{i}\ge0}^{m_{1}+2m_{2}+\cdots+km_{k}=k}\frac{k!}{\vec{m}!}\vec{\nabla}_{\vec{\xi}}^{\left|\vec{m}\right|}\phi(\vec{\xi}):\prod_{i=1}^{k}\left(\frac{\partial_{\vec{d}}^{i}\vec{\xi}(\vec{x}_{0})}{i!}\right)^{m_{i}},
\]
where $\vec{m}=[m_{1},\dots,m_{k}]$, $\vec{m}!=\prod_{i=1}^{k}m_{i}!$,
$\left|\vec{m}\right|=\sum_{i=1}^{k}m_{i}$, $\vec{\nabla}_{\vec{\xi}}^{\left|\vec{m}\right|}$
denotes the derivative tensor of order $\left|\vec{m}\right|$, and
``$:$'' denotes the scalar product of $k$th-order tensors. For
AES-FEM, if (\ref{eq:bilinear-form-gwr-interior}) is used with a
higher-order boundary representation, then the Ciarlet-Raviart condition
is also required when computing $\vec{\nabla}\psi_{i}$. However,
the generalized bilinear form (\ref{eq:bilinear-form-gwr-interior})
uses only $\psi_{i}(\vec{\xi})$ over the master element, so the Ciarlet-Raviart
condition is no longer required. In this case, however, AES-FEM requires
imposing Neumann boundary conditions similar to the techniques in
GFDM instead of FEM. The well-posedness and consistency of such a
hybrid treatment over curved boundaries require a more general analysis,
and we defer it to future work.

\section{\label{sec:convergence} Convergence in $\ell^{p}$ Norm}

In this section, we analyze the convergence of GWR methods in $\ell^{p}$
norm. Given a vector of nodal errors $\delta\vec{u}$, the \emph{convergence
rate in $\ell^{p}$ norm} is $k$th order if $\Vert\delta\vec{u}\Vert_{p}=\sqrt[p]{n}\mathcal{O}(h^{k})$,
where $h$ is some edge length measure of a quasiuniform mesh.

\subsection{Convergence of GWR Methods for Constrained Projection}

Let $\delta\vec{u}_{P-I}=\vec{u}_{P}-\vec{u}_{I}$, where $\vec{u}_{P}$
and $\vec{u}_{I}$ are the projected and interpolated nodal values.
We then obtain the following result regarding the convergence of the
constrained projection.
\begin{thm}
[Convergence of constrained projection]\label{thm:constrained-prjection}Under
the same assumptions as in Theorem~\ref{thm:well-posedness}, if
$u$ is continuously differentiable to $p$th order within each element,
the solution of the constrained projection with a well-posed degree-$p$
GWR method converges \textup{\emph{at }}\textup{$\mathcal{O}(h^{p+1})$}\textup{\emph{
or better in $\ell^{\infty}$ norm, }}i.e., $\left\Vert \delta\vec{u}_{P-I}\right\Vert _{\infty}\leq\mathcal{O}(h^{p+1})$.
\end{thm}
\begin{pf}
Let $\vec{r}=\left\langle u-u_{I},\vec{\Psi}_{\circ;N}\right\rangle _{\Omega}$.
Due to the consistency of GLP basis functions, $\left\Vert \vec{r}\right\Vert _{\infty}\leq\left\Vert \vec{M}\right\Vert _{\infty}\mathcal{O}(h^{p+1})$,
so
\[
\left\Vert \delta\vec{u}_{P-I}\right\Vert _{\infty}\leq\left\Vert \vec{M}^{-1}\right\Vert _{\infty}\left\Vert \vec{r}\right\Vert _{\infty}\leq\left\Vert \vec{M}^{-1}\right\Vert _{\infty}\left\Vert \vec{M}\right\Vert _{\infty}\mathcal{O}(h^{p+1})=\mathcal{O}(h^{p+1}).
\]
\end{pf}
The theorem also applies to other $\ell^{p}$ norms. In Theorem~\ref{thm:constrained-prjection},
we use ``$\leq$'' sign to emphasize that the bound may not be tight
due to possible superconvergence. Note that Theorem~\ref{thm:constrained-prjection}
applies to FEM, GFDM, and AES-FEM. For all of these methods, if $p$
is even, the leading error term in the residual $\vec{r}$ is odd
order, which may cancel out in the integration. In turn, it may lead
to superconvergence. The superconvergence of FEM for the $L^{2}$
projection was considered in \cite{wahlbin1995superconvergence}.

\subsection{\label{subsec:Convergence-of-GWR} Convergence of GWR Methods for
Elliptic BVPs}

Let $\vec{u}_{h}$ denote the nodal solutions of (\ref{eq:linear-system})
corresponding to the nodes in $\Omega_{h}^{\circ}\cup\Gamma_{h,N}^{\circ}$.
Let $\delta\vec{u}_{h-I}=\vec{u}_{h}-\vec{u}_{I}$. We first give
a loose error bound using an argument similar to that of Theorem~\ref{thm:constrained-prjection}.
We state the bound in $\ell^{\infty}$ norm, but the result also holds
in any $\ell^{p}$ norm for $1\leq p\leq\infty$.

\begin{lem}
\label{lem:convergence-GWR-elliptic}Under the same assumptions as
in Theorem~\ref{thm:well-posed-elliptic-bvp}, assuming exact arithmetic,
the solution with a well-posed degree-$p$ GWR method for an elliptic
BVP converges \textup{\emph{at }}\textup{$\mathcal{O}(h^{p-1})$}\textup{\emph{
or better in $\ell^{\infty}$ norm, }}i.e., $\left\Vert \delta\vec{u}_{h-I}\right\Vert _{\infty}\leq\mathcal{O}(h^{p-1})$.
\end{lem}
\begin{pf}
Let $\vec{r}=a(u-u_{I},\vec{\Psi}_{\circ,N})$. Using the same argument
as for Lemma~\ref{lem:deri-approximation}, it is easy to show that
$\left\Vert \mathcal{L}\vec{\Phi}_{i}^{T}\delta\vec{u}_{h-I}\right\Vert _{\infty}=\mathcal{O}(h^{p-1})$
for $\vec{x}_{i}\in\Omega_{h}$. Hence, $\left\Vert \vec{r}\right\Vert _{\infty}\leq\left\Vert \vec{M}\right\Vert _{\infty}\mathcal{O}(h^{p-1})$,
and
\[
\left\Vert \delta\vec{u}_{h-I}\right\Vert _{\infty}\leq\left\Vert \vec{A}^{-1}\right\Vert _{\infty}\left\Vert \vec{r}\right\Vert _{\infty}\leq\left\Vert \vec{A}^{-1}\right\Vert _{\infty}\left\Vert \vec{M}\right\Vert _{\infty}\mathcal{O}(h^{p-1})=\mathcal{O}(h^{p-1}).
\]
\end{pf}
\begin{rem}
Lemma~\ref{lem:convergence-GWR-elliptic} assumes exact arithmetic.
This is important because in the presence of approximation or rounding
errors, the solutions may not converge due to \emph{ill-conditioning};
see e.g. \cite[p. 45]{leveque2007finite} for finite differences and
\cite[p. 222]{ern2013theory} for finite elements. 
\end{rem}
For FEM, Lemma~\ref{lem:convergence-GWR-elliptic} underestimates
the convergence rate by two orders, compared to the well-known $\mathcal{O}(h^{p+1})$
error bounds in $L^{2}$ norm due to the Aubin-Nitsche duality argument;
see e.g. \cite{strang1973analysis}. The derivation of error bounds
for FEM in $\ell^{p}$ norm is beyond the scope of this work. For
GFDM and AES-FEM, however, the Aubin-Nitsche duality argument does
not apply, due to the non-conformity. For odd-degree-$p$ GFDM and
AES-FEM, Lemma~\ref{lem:convergence-GWR-elliptic} is tight. However,
for even-degree-$p$, similar to constrained projections, the leading
error term in the residual $\vec{r}$ is odd order, which may cancel
out in the integration. For GFDM, this error cancellation occurs with
symmetric stencils, analogous to centered differences, leading to
$\mathcal{O}(h^{p})$ convergence rate. With AES-FEM, however, the
error cancellation is primarily due to the numerical integration,
for example, see Figure \ref{fig:convergence-solution}.
\begin{thm}
\label{thm:AES-FEM-Convergence} Under the same assumptions as in
Theorem~\ref{thm:well-posed-elliptic-bvp}, assuming exact arithmetic,
the solution of a well-posed even-degree-$p$ AES-FEM for a coercive
elliptic BVP with Dirichlet boundary conditions converges \textup{\emph{at
$\mathcal{O}(h^{p})$ in $\ell^{\infty}$ norm, }}if $p=2$ or if
$p\ge4$ and the local support is nearly symmetric.
\end{thm}
\begin{pf}
Let $u_{I,i}$ denote the local interpolation at a node $\vec{x}_{i}\in\Omega_{h}^{\circ}$.
Let $\delta u_{i}=u_{I,i}-u$, which is a smooth function over $\Omega_{i}$.
Apply $\mathcal{L}$ to $\delta u_{i}$. We note that 
\[
\mathcal{L}\delta u_{i}(\vec{x}_{i}+\vec{h})=C\left(\text{sign}\left(\vec{h}\cdot\hat{\vec{h}}\right)\left\Vert \vec{h}\right\Vert \right)^{p-1}+\text{\ensuremath{\mathcal{O}}}(h^{p}),
\]
where $C$ is proportional to $\left\Vert \partial_{\hat{\vec{h}}}^{p+1}u(\vec{x})\right\Vert _{\infty}$.
If $p=2$, because $\psi_{i}$ is the hat function, the line integral
of $\left(\text{sign}\left(\vec{h}\cdot\hat{\vec{h}}\right)\left\Vert \vec{h}\right\Vert \right)^{p-1}\psi_{i}(\vec{x}_{i}+\vec{h})$
cancels out exactly, so $\left\Vert \vec{r}\right\Vert _{\infty}\leq\left\Vert \vec{M}\right\Vert _{\infty}\mathcal{O}(h^{2})$
in Lemma~\ref{lem:convergence-GWR-elliptic}. For even $p\geq4$
, the line integral cancels out if the local support is (nearly) symmetric
about $\vec{x}_{i}$.
\end{pf}
In \cite{conley_2016}, quadratic AES-FEM was shown to converge at
second order, but the proof did not explicitly state the error cancellation.
Theorem~\ref{thm:AES-FEM-Convergence} indicates that AES-FEM may
not enjoy the full $\mathcal{O}(h^{p})$ superconvergence for $p\ge4$
on highly irregular meshes, but in practice we observe it to be close
to $\mathcal{O}(h^{p})$ on quasiuniform meshes, as we will demonstrate
in Section~\ref{subsec:Convergence-of-High-Order-AES-FEM}.

\section{\label{sec:numerical-results} Numerical Results}

In this section, we present some numerical results to verify the theoretical
analysis in this work.

\subsection{\label{subsec:Element-Shape-Independence} Comparison of Mesh Dependency}

\begin{figure}
\begin{minipage}[t]{0.48\textwidth}%
\begin{center}
\includegraphics[width=1\columnwidth]{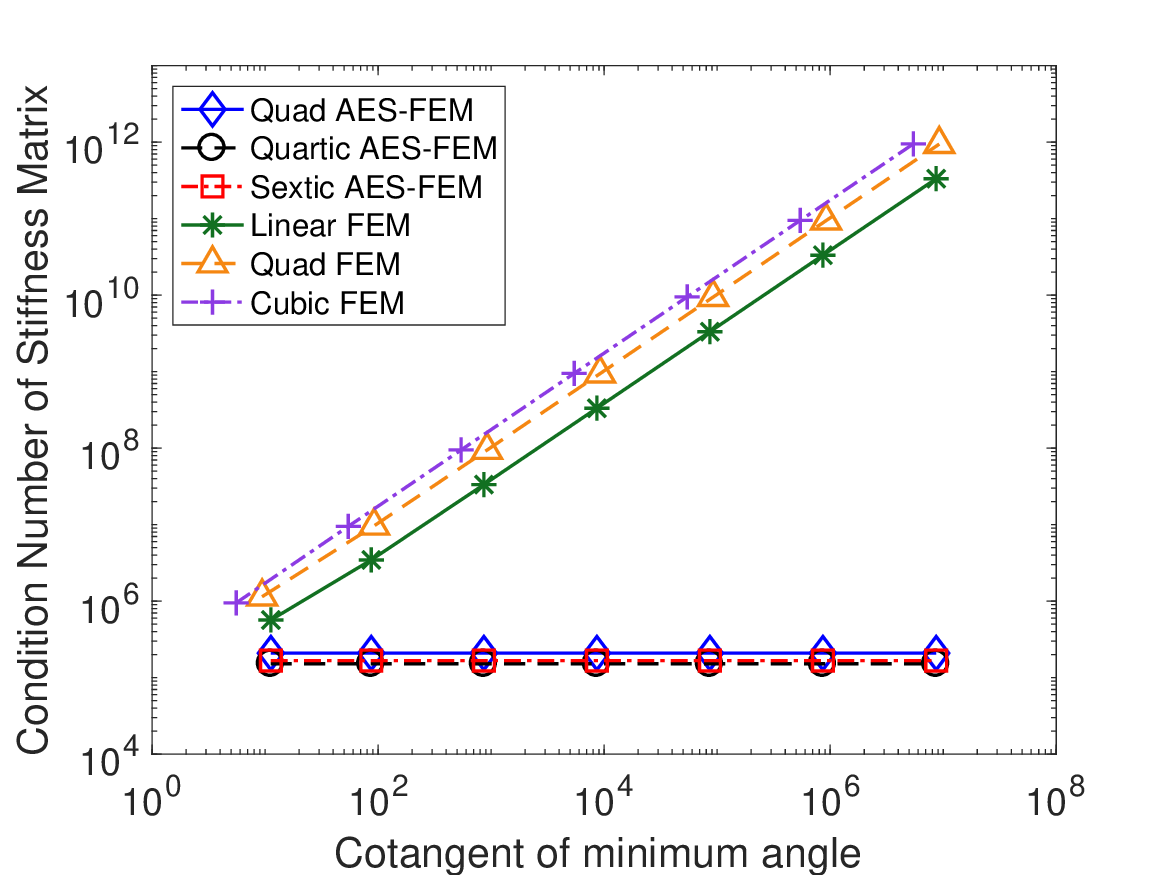}
\par\end{center}
\begin{flushleft}
\caption{\label{fig:conds_bad_mesh}Dependence of the condition numbers of
the stiffness matrices of FEM and AES-FEM on the worse angles.}
\par\end{flushleft}%
\end{minipage}\hfill%
\begin{minipage}[t]{0.48\textwidth}%
\begin{center}
\includegraphics[width=1\columnwidth]{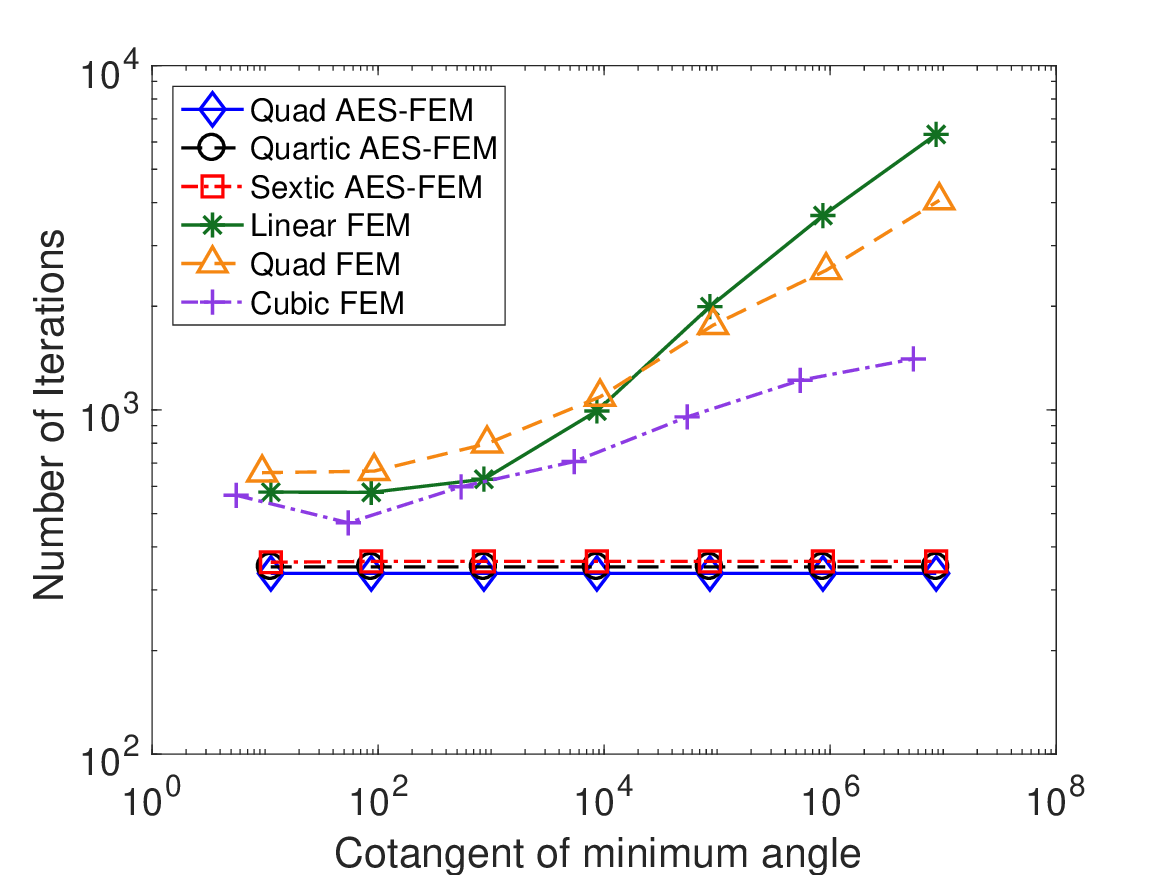}
\par\end{center}
\begin{flushleft}
\caption{\label{fig:iters_bad_mesh_2d}Relationship of the number of iterations
of the preconditioned iterative solvers and the worst angles.}
\par\end{flushleft}%
\end{minipage}
\end{figure}

We first compare the mesh dependency of FEM to that of AES-FEM. As
shown in Section~\ref{subsec:Mesh-Dependency}, the well-posedness
of FEM depends on well-shapedness of the meshes, while AES-FEM is
independent of element shapes, regardless of the degree of the GLP
basis functions. To demonstrate this, we solved the Poisson equation
in 2D
\begin{equation}
-\Delta u=f\label{eq:Poisson-eq}
\end{equation}
with Dirichlet boundary conditions over $[-1,1]^{2}$. We obtained
$f$ and $u_{D}$ from the exact solution $u=\cos(\pi x)\cos(\pi y)$
and solved the equation using our own implementations of quadratic,
quartic, and sextic AES-FEM, along with linear, quadratic, and cubic
FEM. We generated a series of meshes on a unit square with progressively
worse element quality, which we obtain by distorting a good-quality
mesh. For AES-FEM and linear FEM, we used a mesh with 130,288 elements
and 65,655 nodes and distorted four elements by moving one vertex
of each of these elements incrementally towards its opposite edge.
For quadratic FEM, the mesh had 32,292 elements and 65,093 nodes,
and a single element was distorted by moving one vertex and its adjacent
mid-edge nodes incrementally towards its opposite edge. For cubic
FEM, the mesh had 32,292 elements and 146,077 nodes, and also a single
element was distorted. Figure~\ref{fig:conds_bad_mesh} shows the
condition numbers of the stiffness matrices of FEM and AES-FEM.

In practice, the condition number may affect the efficiency of iterative
solvers. Figure~\ref{fig:iters_bad_mesh_2d} shows the numbers of
iterations required to solve the linear systems to a relative tolerance
of $10^{-8}$ using GMRES for AES-FEM and CG for FEM, both with Gauss-Seidel
preconditioners. It can be seen that the condition numbers of FEM
increased inversely proportional to the minimum angle, and the number
of iterations of CG grew correspondingly. In contrast, the condition
numbers and the number of iterations for AES-FEM remained constant.
We observed similar behavior for 3D AES-FEM, which we omit from the
paper.

\subsection{\label{subsec:Convergence-of-High-Order-AES-FEM} Convergence of
High-Order AES-FEM}

\begin{figure}
\begin{minipage}[t]{0.33\textwidth}%
\begin{center}
\includegraphics[width=1\textwidth]{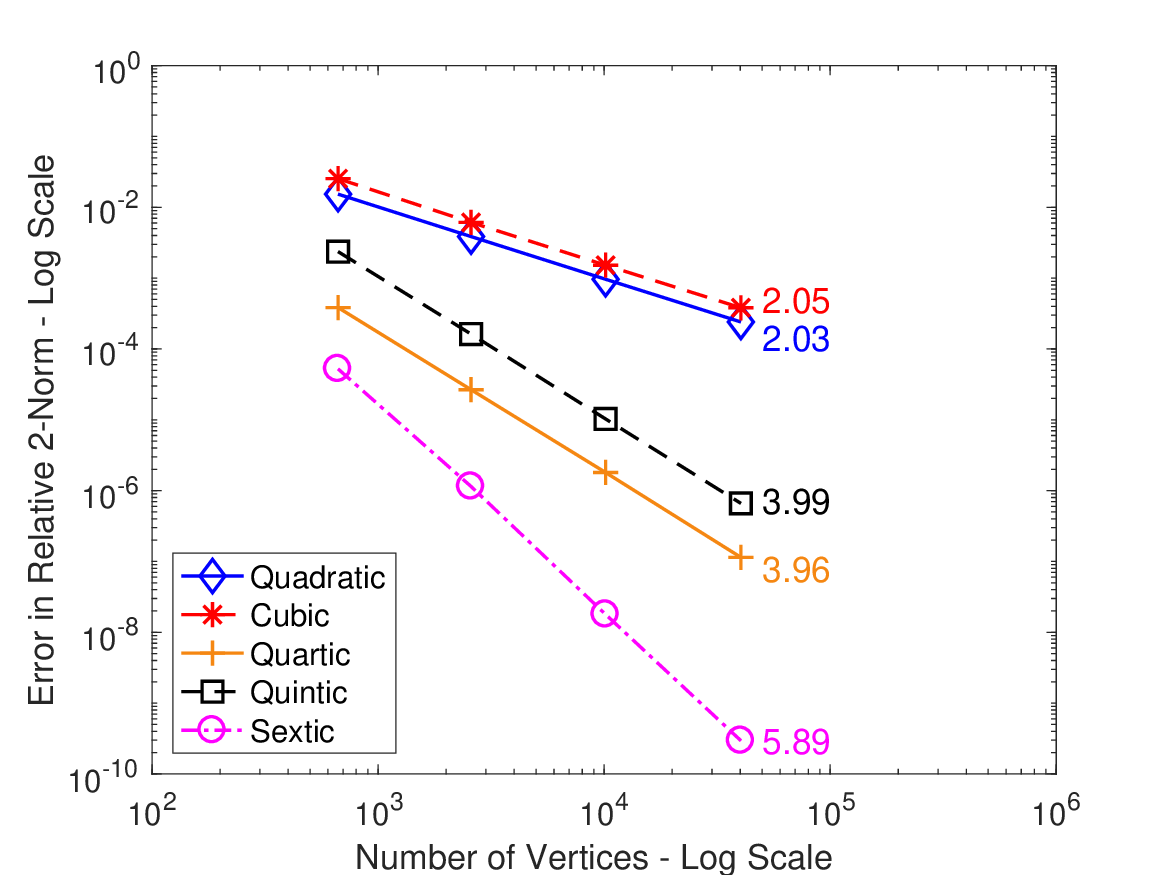}\\
{\footnotesize{}(a) Poisson equation }\\
{\footnotesize{} with Dirichlet boundaries.}
\par\end{center}%
\end{minipage}%
\begin{minipage}[t]{0.33\textwidth}%
\begin{center}
\includegraphics[width=1\textwidth]{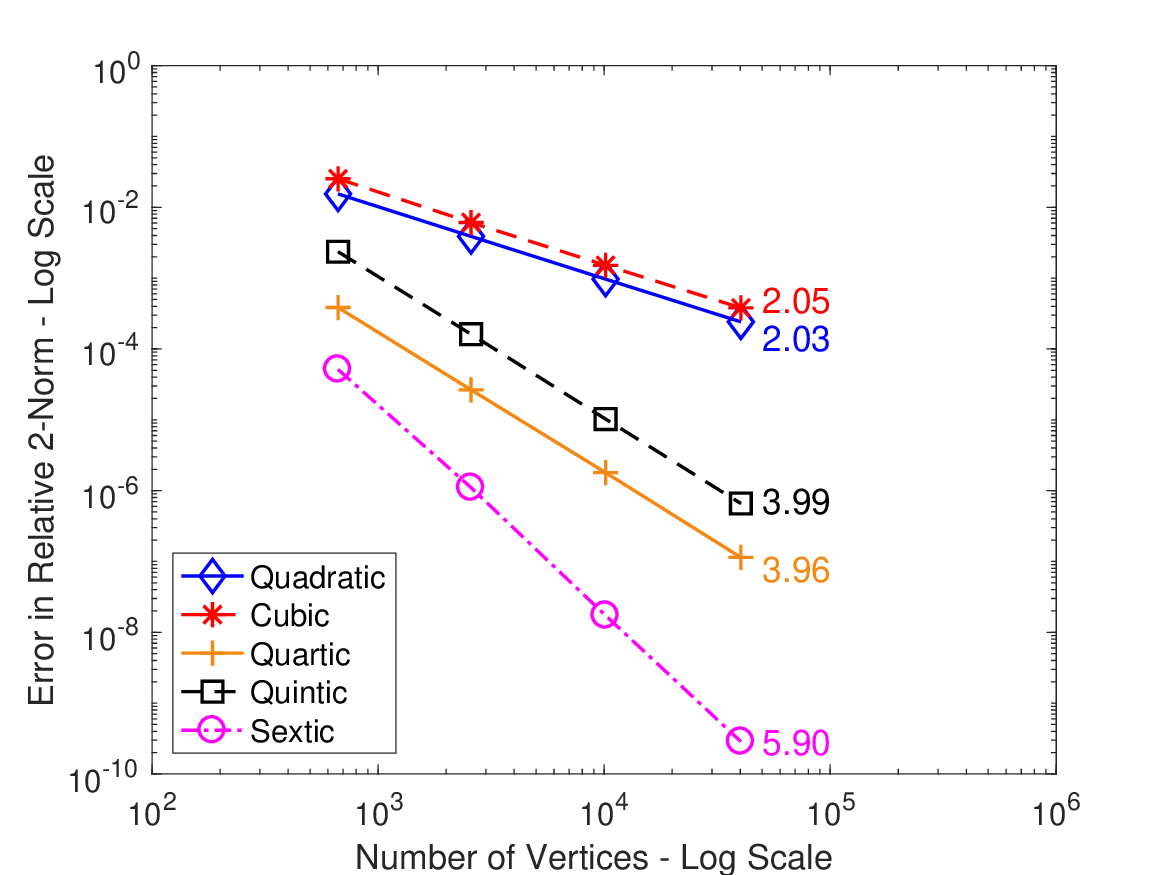}\\
{\footnotesize{}(b) Advection-diffusion equation }\\
{\footnotesize{} with Dirichlet boundaries.}
\par\end{center}%
\end{minipage}%
\begin{minipage}[t]{0.33\textwidth}%
\begin{center}
\includegraphics[width=1\textwidth]{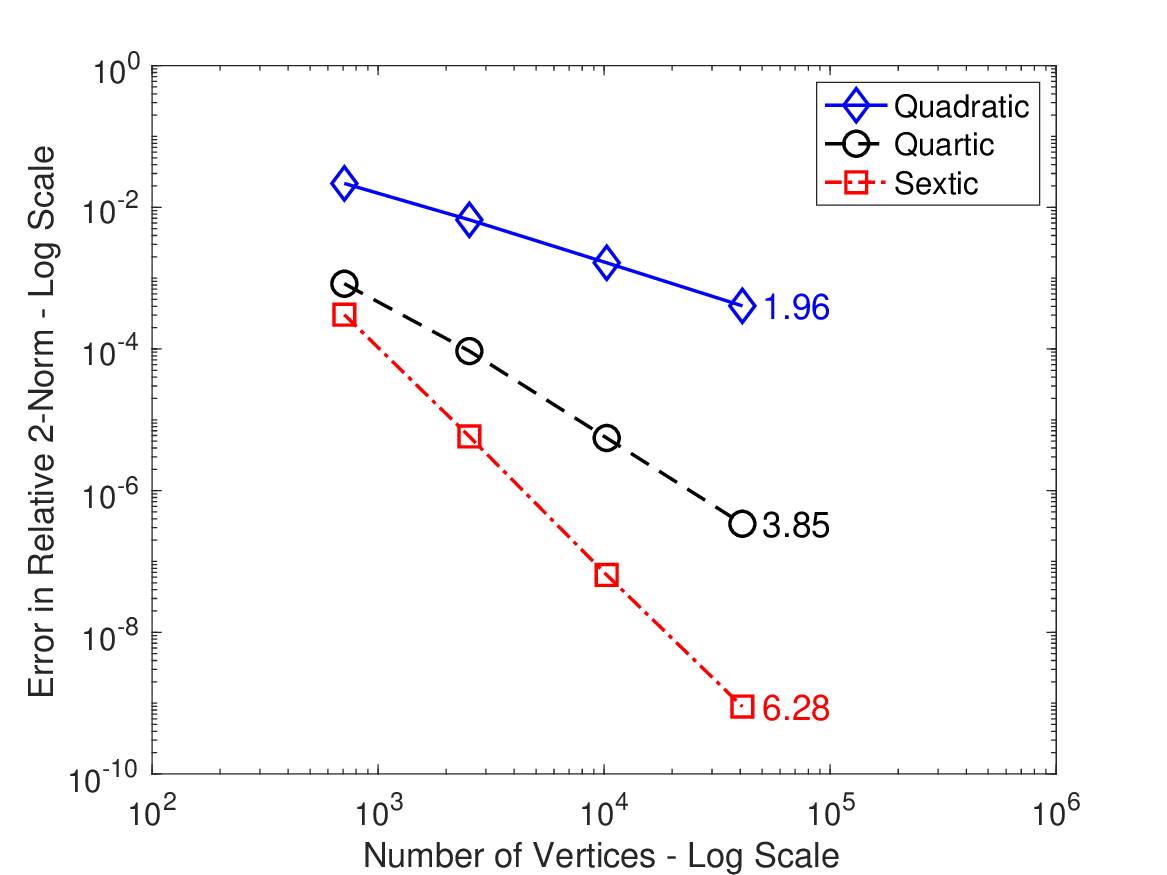}\\
{\footnotesize{}(c) Advection-diffusion equation}\\
{\footnotesize{} with Neumann boundaries.}
\par\end{center}%
\end{minipage}
\raggedright{}\caption{\label{fig:convergence-solution} The errors in the function values
from AES-FEM for 2D Poisson and advection-diffusion equations. The
number to the right of each curve indicates the average convergence
rate.}
\end{figure}

Next, we verify the convergence analysis in Section~\ref{subsec:Convergence-of-GWR},
especially that of AES-FEM. To this end, we used AES-FEM of degrees
2 to 6 to solve the equation
\begin{equation}
-\Delta u+\vec{\nu}\cdot\vec{\nabla}u=f,\label{eq:convection-diffusion}
\end{equation}
with unstructured meshes over $[-1,1]^{2}$, where $f=\sin\pi x\sin\pi y$.
Figure~\ref{fig:convergence-solution}(a) shows the convergence rates
in the relative $\ell^{2}$ norm for the Poisson equation, for which
$\vec{\nu}=\vec{0}$. Figure~\ref{fig:convergence-solution}(b) and
(c) show the convergence rates for the advection-diffusion equation
with Dirichlet and Neumann boundary conditions, respectively, where
$\vec{\nu}=[x,-y]$. The number to the right of each convergence curve
shows the average convergence rate under mesh refinement. Note that
AES-FEM with even-degree basis functions converged at about $p$th
order whereas with odd-degree basis functions, AES-FEM converged at
$(p-1)$st order. For example, with quadratic and cubic basis functions,
the convergence rate is approximately second order. This difference
in convergence rates is due to error cancelation in the numerical
integration, as discussed in Section \ref{subsec:Convergence-of-GWR}.

\subsection{\label{subsec:Comparison-with-FEM}Comparison of FEM, GFDM, and AES-FEM}

\begin{figure}
\begin{centering}
\begin{minipage}[t]{0.33\textwidth}%
\begin{center}
\includegraphics[width=0.7\columnwidth]{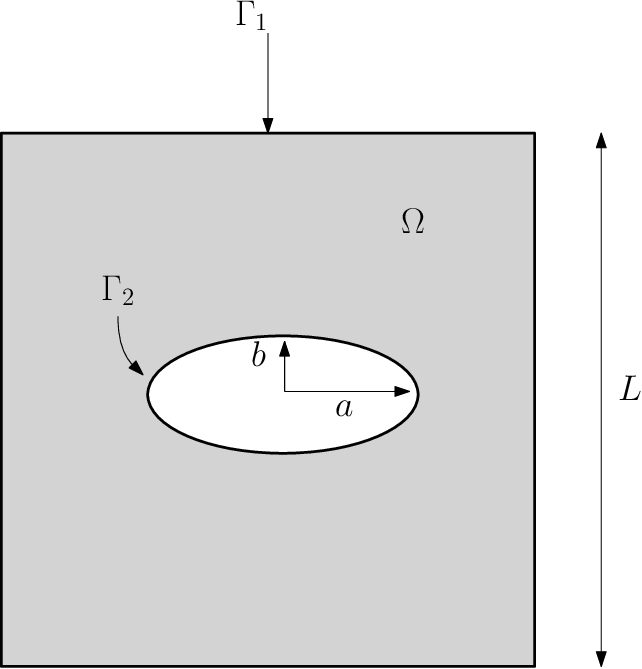}\\
{\footnotesize{}(a) Test domain.}
\par\end{center}%
\end{minipage}%
\begin{minipage}[t]{0.33\textwidth}%
\begin{center}
\includegraphics[width=1\columnwidth]{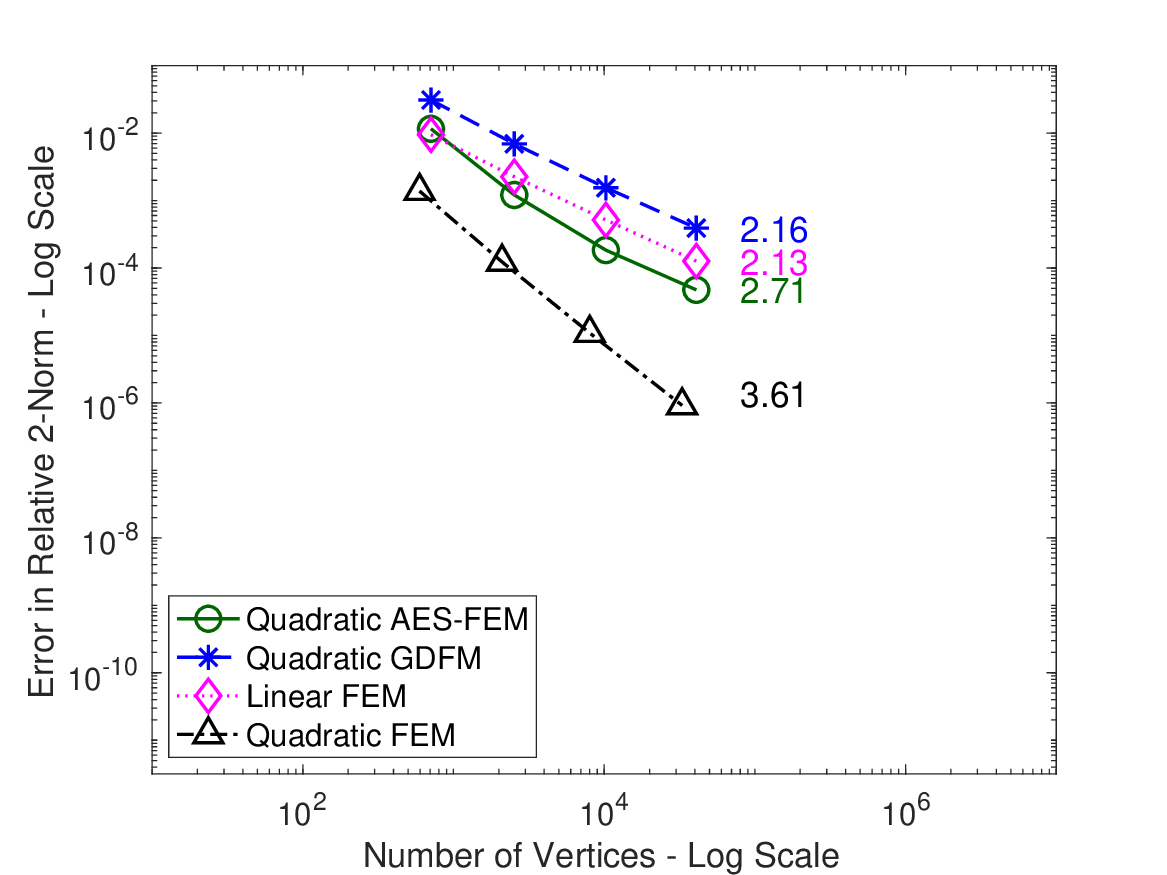}\\
{\footnotesize{}(b) Low-order methods.}
\par\end{center}%
\end{minipage}%
\begin{minipage}[t]{0.33\textwidth}%
\begin{center}
\includegraphics[width=1\columnwidth]{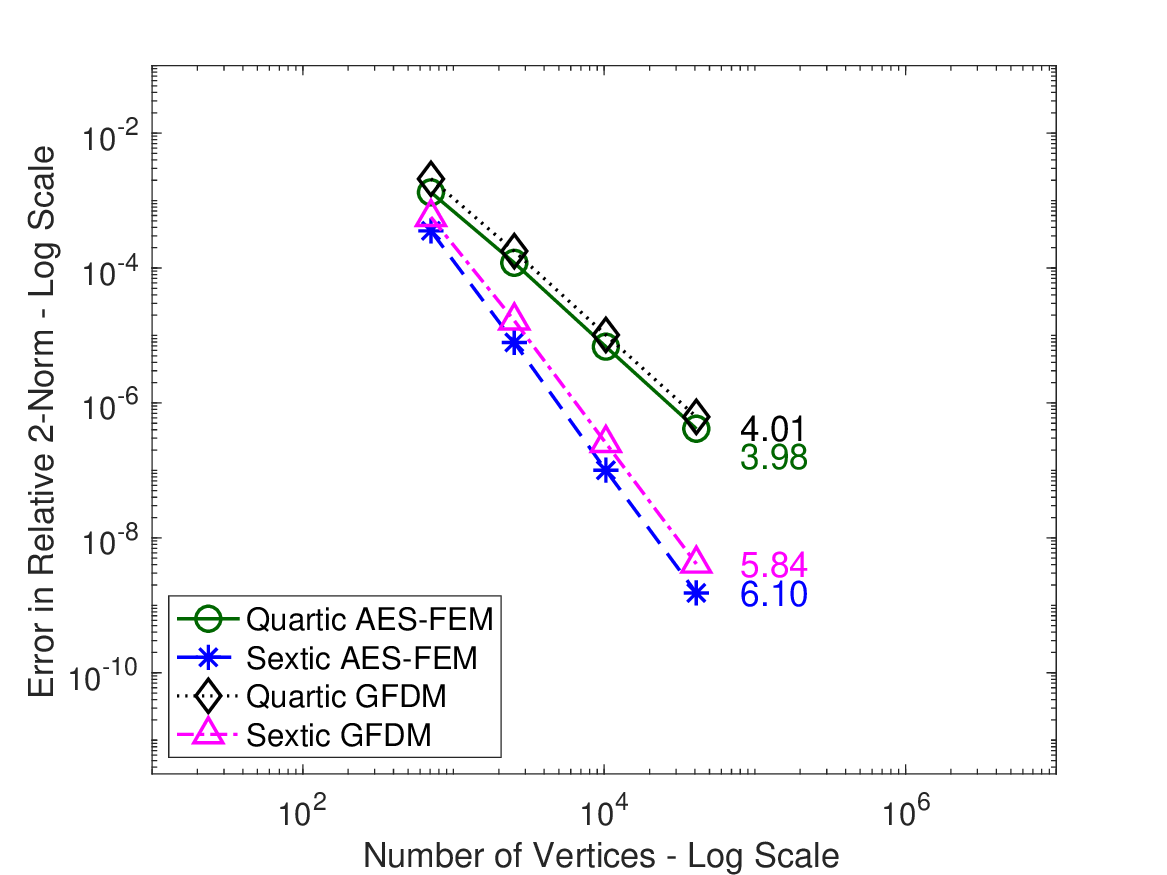}\\
{\footnotesize{}(c) High-order methods.}
\par\end{center}%
\end{minipage}
\par\end{centering}
\caption{\label{fig:comparison-of-FEM}Comparison of FEM, GFDM, and AES-FEM
for the Poisson equation.}
\end{figure}

Finally, we compare the accuracy of FEM, GFDM, and FEM. We use the
Poisson equation in this comparison. In particular, we solved a 2D
Poisson equation (\ref{eq:Poisson-eq}) over the square $[-1,1]^{2}$
with an elliptical hole of semi-axes $0.5$ and $0.2$ in the middle.
The domain, as illustrated in Figure~\ref{fig:comparison-of-FEM}(a),
has nonuniform curvature along the inner boundary and has corners
along the outer boundary. We applied Neumann boundary conditions to
the outer boundary and applied Dirichlet boundary conditions to the
inner boundary. We obtained the source term $f$ and the boundary
conditions by differentiating the following analytic function
\begin{equation}
u=\sin\left(\pi x\right)\sin\left(\pi y\right).\label{eq:Poisson_sol_2D}
\end{equation}

First, let us focus on comparing linear and quadratic FEM, GFDM, and
AES-FEM. For FEM, we used the Partial Differential Equation (PDE)
Toolbox in MATLAB R2018a \cite{matlab_pde2018}, which supports linear
and quadratic elements. Hence, we focus on comparing linear and quadratic
FEM with quadratic AES-FEM and GFDM. We refer to all these methods
as ``low-order''. We generated the meshes directly using PDE Toolbox
and used the built-in solvers in MATLAB with default tolerances. The
number of nodes ranged between 709 and 40,872. Figure~\ref{fig:comparison-of-FEM}(b)
compares these lower-order methods, where the number to the right
of each convergence curve indicates the average convergence rates.
Quadratic AES-FEM had slightly better accuracy than linear FEM on
finer meshes. GFDM and AES-FEM have identical sparsity patterns, but
GFDM had significantly larger errors than AES-FEM.

Next, we compare quartic and sextic AES-FEM and GFDM. We used the
same meshes as for quadratic AES-FEM. Figure~\ref{fig:comparison-of-FEM}(c)
compares AES-FEM and GFDM for the problem. For GFDM, we applied Neumann
conditions by averaging the one-sided derivatives on both left- and
right-hand sides at corners. It can be seen that AES-FEM slightly
outperformed GFDM in all the cases.

It should be noted that the linear systems for quadratic AES-FEM and
linear FEM have nearly identical sparsity patterns. Additionally,
the linear systems of quartic AES-FEM has only slightly more nonzeros
than that of quadratic FEM, but it is significantly more accurate.
Thus, when comparing matrices with similar numbers of DOFs and similar
numbers of nonzeros, AES-FEM is often more accurate than FEM. As a
result, AES-FEM sometimes requires less computational time than FEM
\cite[Chapter 6]{conley_2016,conley2016overcoming}.

\section{\label{sec:Conclusions} Conclusions and Discussions}

In this paper, we introduced the framework of \emph{generalized weighted
residual formulations} (\emph{GWR}), which unifies generalized finite
differences, Lagrange finite elements, and adaptive extended stencil
FEM (AES-FEM). Under this framework, we presented a unified analysis
of the \emph{well-posedness} of these methods, which depend on the
\emph{quasiuniformity} and the consistency and stability of the trial
functions and test functions. While the stability of the Lagrange
basis functions in FEM depends on the well-shapedness of the elements,
the GLP basis functions of GFDM and AES-FEM depend on the selections
of the stencils, which can be adapted locally, due to the least-squares
nature of the GLP basis functions. In addition, high-order AES-FEM
requires only first-order meshes for its implementation, so its implementation
is simpler than high-order finite elements. However, Lagrange FEM
can achieve $\mathcal{O}(h^{p+1})$ convergence rate in $L^{2}$ norm.
In contrast, GFDM and AES-FEM significantly simplify mesh generation,
but it comes at the cost of a lower-order convergence rate, which is
$\mathcal{O}(h^{p-1})$ with odd-degree GLP basis functions. However,
with even-degree basis functions, GFDM can achieve $\mathcal{O}(h^{p})$
convergence rate with nearly symmetric stencils, whereas AES-FEM can
achieve $\mathcal{O}(h^{p})$ convergence with nearly symmetric local
support. We presented numerical results to verify our theoretical
analysis, and we showed that AES-FEM in general outperforms GFDM in
terms of accuracy. For Neumann boundary conditions, we only considered
polygonal domains. For curved geometries, we can overcome mesh quality
dependency, namely the Ciarlet-Raviart condition, in AES-FEM by using
techniques similar to GFDM, but a rigorous proof of well-posedness
and consistency is challenging, which we plan to report elsewhere.
One direction of future work is to investigate the development of
new hybrid methods by mixing FEM and AES-FEM for different orders
of terms.

\section*{Acknowledgments}

This work was supported in part by DoD-ARO under contract \#W911NF0910306
and also in part under the Scientific Discovery through Advanced Computing
(SciDAC) program in the US Department of Energy Office of Science,
Office of Advanced Scientific Computing Research  through subcontract
\#462974 with Los Alamos National Laboratory and  under a subcontract
with Argonne National Laboratory under Contract DE-AC02-06CH11357.
The first author acknowledges the support of the Kenny Fund Fellowship
of Saint Peter's University. Results were obtained using the high-performance
LI-RED computing system at the Institute for Advanced Computational
Science at Stony Brook University, which was obtained through the
Empire State Development grant NYS \#28451.
\bibliography{refs/aes_fem_high}

\appendix

\section{\label{sec:Computation-of-GLP}Computation of GLP Basis Functions}

Given a node $\vec{x}_{0}$, let $\vec{v}$ denote the local coordinate
system centered at it. Let $\vec{\mathcal{P}}_{k}^{(p)}(\vec{x})$
denote the set of all $d$-dimensional monomials of degree $p$ and
lower; for example, $\vec{\mathcal{P}}_{2}^{(2)}(\vec{x})=[1,x,y,x^{2},xy,y^{2}]^{T}$.
Let $\vec{D}_{k}^{(p)}$ be a diagonal matrix consisting of the fractional
factorial part of the coefficients in the Taylor series corresponding
to $\vec{\mathcal{P}}_{k}^{(p)}$; for example, $\vec{D}_{2}^{(2)}=\text{diag}\left(1,1,1,\nicefrac{1}{2},1,\nicefrac{1}{2}\right)$.
Let $\vec{c}$ be a vector containing the partial derivative of $f$
evaluated at $\vec{x}_{0}$; for example, $\vec{c}=\left[f,f_{x},f_{y},f_{xx},f_{xy},f_{yy}\right]^{T}\mid_{\vec{x}=\vec{x}_{0}}$.
Then, we may write the truncated Taylor series of a smooth function
$f$ as
\begin{equation}
f(\vec{x})\approx\vec{c}^{T}\vec{D}_{k}^{(p)}\vec{\mathcal{P}_{k}}^{(p)}\left(\vec{x}\right).\label{eq:TaylorSeries}
\end{equation}
Suppose there are $n$ coefficients in $\vec{c}$, and the stencil
about the point $\vec{x}_{0}$ contains $m$ points, including the
point $\vec{x}_{0}$. To obtain the $j$th basis function $\phi_{j}$,
let $f(\vec{x}_{i})=\delta_{ij}$, the Kronecker delta function. Therefore,
we obtain an $m\times n$ least squares problem 
\begin{equation}
\vec{V}\vec{c}_{j}\approx\vec{e}_{j},\label{eq:least-squares}
\end{equation}
where $\vec{e}_{j}$ denotes the $j$th column of the $m\times m$
identity matrix, and $\vec{V}$ is the generalized Vandermonde matrix.
Eq.~(\ref{eq:least-squares}) may potentially be ill-conditioned
and potentially rank deficient, even if $m\ge n$. We solve (\ref{eq:least-squares})
by minimizing a weighted norm (or semi-norm) 
\begin{equation}
\min_{\vec{c}}\left\Vert \vec{V}\vec{c}_{j}-\vec{e}_{j}\right\Vert _{\vec{W}}\equiv\min_{\vec{c}}\left\Vert \vec{W}\left(\vec{V}\vec{c}_{j}-\vec{e}_{j}\right)\right\Vert _{2},
\end{equation}
where $\vec{W}$ is an $m\times m$ diagonal weighting matrix, and
it is a constant for a given node. In general, heavier weights are
assigned to nodes that are closer to $\vec{x}_{0}$; for example,
\begin{equation}
w_{i}=\left(\frac{\left\Vert \vec{v}_{i}\right\Vert }{h}+\epsilon\right)^{-p/2},\label{eq:weighting_scheme_AES}
\end{equation}
where $\epsilon$ is a small number, such as $\epsilon=0.01$, for
avoiding division by zero.

The matrix $\vec{WV}$ can be poorly scaled. We address this by right-multiplying
a diagonal matrix $\vec{S}$. Let $\vec{a}_{j}$ denote the $j\text{th}$
column of an arbitrary matrix $\vec{W}\vec{V}$. A typical choice
for the $i\text{th}$ entry of $\vec{S}$ is either $1/\left\Vert \vec{a}_{i}\right\Vert _{2}$
or $1/\left\Vert \vec{a}_{i}\right\Vert _{\infty}$. This is known
as \emph{column}\textit{ equilibration}\textit{\emph{ \cite{van1969condition}.
Note that a row equilibration or a general matrix equilibration cannot
be used, since it would undermine the weighting scheme $\vec{W}$.}}
After weighting and scaling, the least-squares problem becomes 
\begin{equation}
\min_{\vec{d}}\left\Vert \tilde{\vec{V}}\vec{d}-\vec{W}\vec{e}_{j}\right\Vert _{2},\quad\text{where }\tilde{\vec{V}}\equiv\vec{W}\vec{V}\vec{S}\text{ and }\vec{d}\equiv\vec{S}^{-1}\vec{c}_{j}.\label{eq:WLS-1}
\end{equation}
We solve the problem using the truncated QR factorization with column
pivoting, where the pivoting scheme is customized to preserve low-degree
terms. The solution of the least squares problem is $\vec{c}_{j}=\vec{S}\tilde{\vec{V}}^{+}\vec{W}\vec{e}_{j}$.
The complete set of basis functions is then given by 
\begin{equation}
\vec{\Phi}=\left(\vec{S}\tilde{\vec{V}}^{+}\vec{W}\right)^{T}\vec{D}\vec{\mathcal{P}}.\label{eq:GLPBF_def}
\end{equation}

\section{\label{sec:Selection-of-Stencils} Selection of Stencils}

\begin{figure}
\begin{minipage}[t]{0.48\textwidth}%
\begin{center}
\includegraphics[width=1\columnwidth]{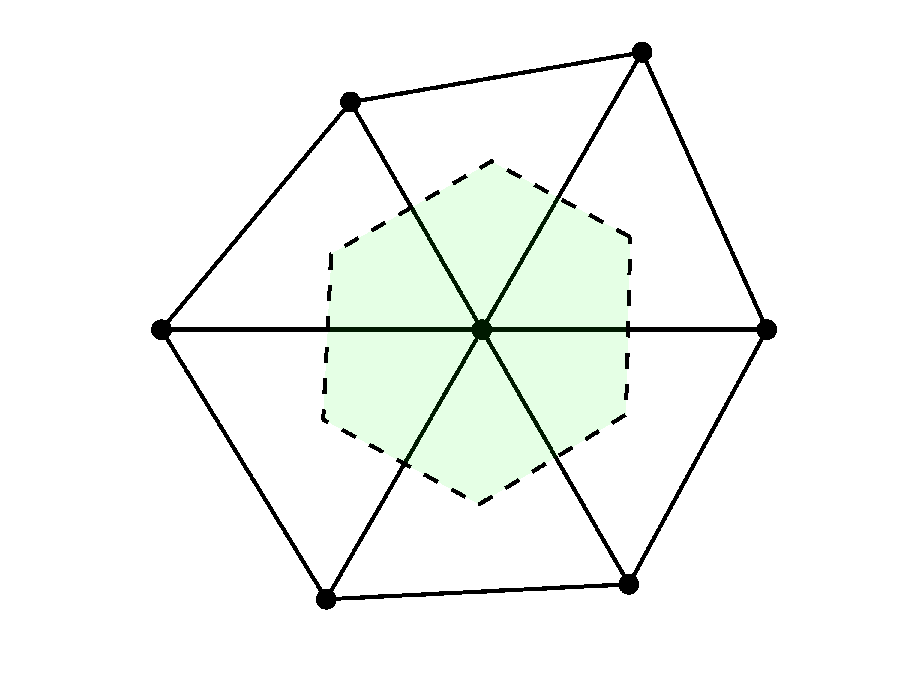}
\par\end{center}
\begin{flushleft}
\caption{\label{fig:control_volume}Example of 1-ring neighbor elements with
the control volume of the center node in light green.}
\par\end{flushleft}%
\end{minipage}\hfill%
\begin{minipage}[t]{0.48\textwidth}%
\begin{center}
\includegraphics[width=0.8\columnwidth]{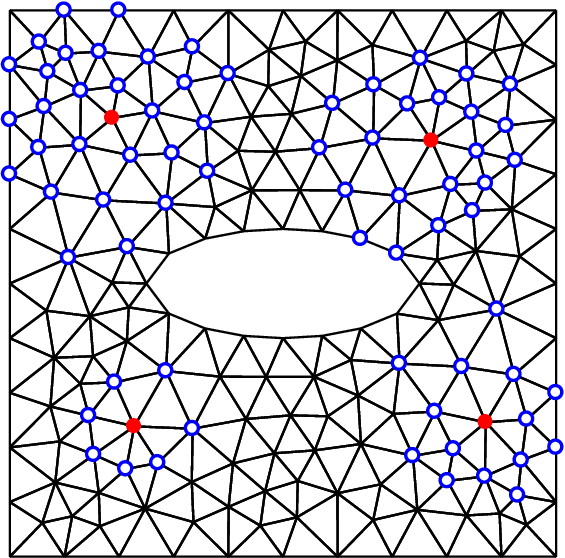}
\par\end{center}
\begin{flushleft}
\caption{\label{fig:stencils_2d}Examples of 2D stencils (blue circles) of
four nodes (red dots) with neighborhood sizes of $2\sfrac{1}{2}$-ring,
$2$-ring, $1\sfrac{1}{2}$-ring, $1$-ring neighborhoods, clockwise
from top left.}
\par\end{flushleft}%
\end{minipage}
\end{figure}

To achieve high-order accuracy, a critical question is the selection
of the stencils at each node for the construction of the GLP basis
functions. We utilize meshes for speedy construction of the stencils.
Given a simplicial mesh (i.e., a triangle mesh in 2D or a tetrahedral
mesh in 3D), the \textit{1-ring neighbor elements} of a node are defined
to be the elements incident on the node. See Figure~\ref{fig:control_volume}
for an example of the 1-ring neighborhood elements and the control
volume of a node. The \textit{1-ring neighborhood} of a node contains
the nodes of its 1-ring neighbor elements \cite{Jiao2008}. For any
integer $k\geq1$, we define the $(k+1)$-\textit{ring neighborhood}
as the nodes in the $k$-ring neighborhood plus their 1-ring neighborhoods. 

\begin{figure}
\begin{centering}
\includegraphics[scale=0.6]{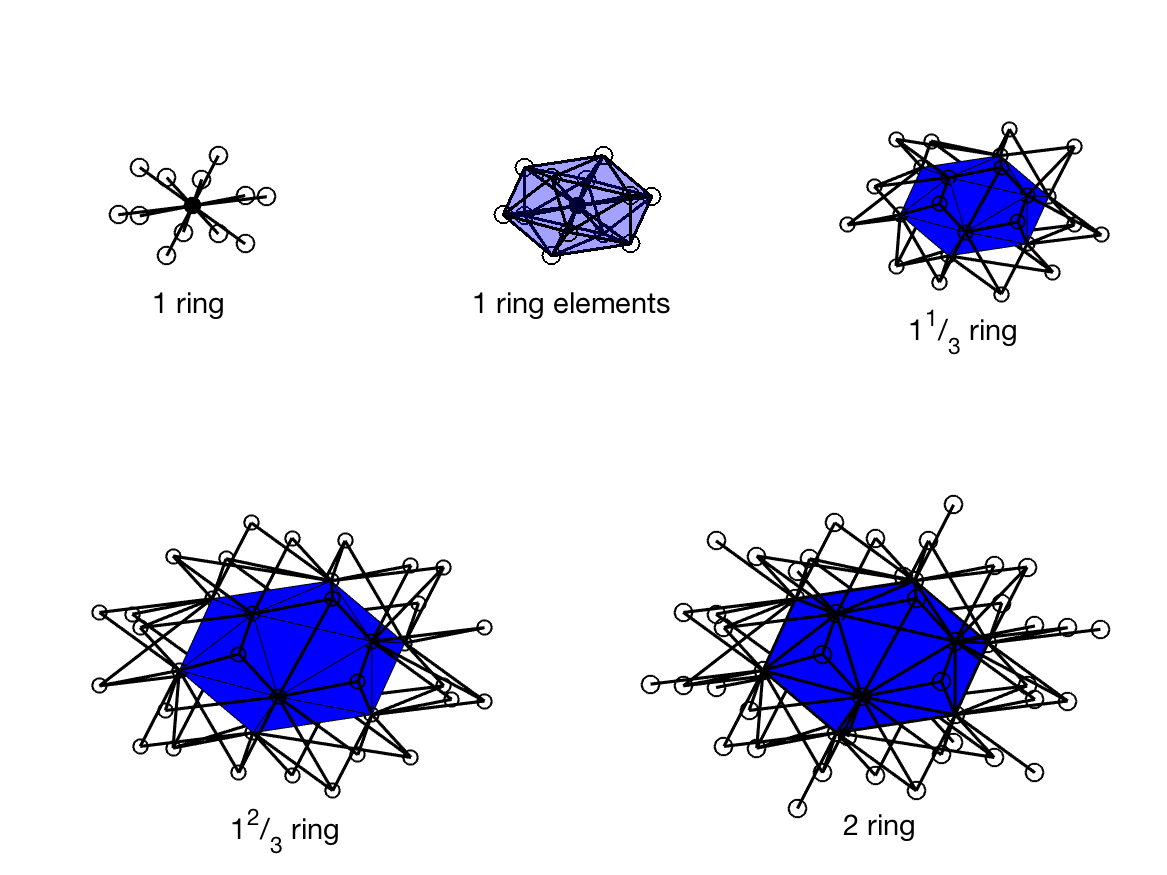}
\par\end{centering}
\caption{Examples of 3D stencils with $1$-ring neighborhood, $1$-ring neighbor
elements, $1\sfrac{1}{3}$-ring neighborhood, $1\sfrac{2}{3}$-ring
neighborhood and $2$-ring neighborhood of the center node (in solid
black).\label{fig:stencil-3D}}
\end{figure}

The 1-ring neighborhood of a node may supply a sufficient number of
nodes for constructing quadratic GLP basis functions. However, 2-\textit{\emph{
and 3-rings are often too large for cubic and quartic constructions.
We refine the}} granularity of the stencils by using fractional rings.
In 2D we use half-rings, as defined in \cite{Jiao2008}. For an integer
$k\geq1$, the $(k+\sfrac{1}{2})$\textit{-ring neighborhood} is the
$k$-ring neighborhood together with the nodes of all the faces that
share an edge with the $k$-ring neighborhood. See Figure~\ref{fig:stencils_2d}
for an illustration. For 3D, we use $\sfrac{1}{3}$- and $\sfrac{2}{3}$-rings,
as defined in \cite{conley_2016}. For any integer $k\geq1$, the
$\left(k+\sfrac{1}{3}\right)$-\textit{ring neighborhood} contains
the $k$-ring neighborhood together with the nodes of all elements
that share a face with the $k$-ring neighborhood. The $(k+\sfrac{2}{3})$-\textit{ring
neighborhood} contains the $k$-ring neighborhood together with the
nodes of all faces that share an edge with the $k$-ring neighborhood.
See Figure~\ref{fig:stencil-3D} for an illustration of rings, one-third
rings and two-third rings in 3D.

In practice, for degree-$p$ basis functions in $d$-dimensional space,
we typically choose the ring size $\left(p+1\right)/d$. This offers
a good balance of accuracy, stability, and efficiency for basis functions
up to degree 7. To illustrate, Table~\ref{tab:m_vs_n} compares the
average number of nodes in a given sized ring to the number of unknowns
for a given degree in (\ref{eq:least-squares}) on an example mesh.
It can be seen that the $\left(p+1\right)/d$-ring size offers approximately
$1.3$ to $2$ times the number of the coefficients on average.
If a particular neighborhood does not provide enough points, especially
for nodes near boundaries, we further expand the stencil to a larger
ring. The construction of the neighborhood requires an efficient mesh
data structure, such as the Array-based Half-Facet (AHF) data structure
\cite{Dyedov2014}.

\begin{table}
\caption{Comparison of the average number of nodes per ring versus the number
of coefficients for 2D (left) and 3D (right) Taylor polynomials.\label{tab:m_vs_n}}
\begin{minipage}[t]{0.49\textwidth}%
\begin{center}
\begin{tabular}{c|c||c|c}
\hline 
Degree & \#Coeffs. & Ring & \#Nodes\tabularnewline
\hline 
\hline 
2 & 6 & $1\sfrac{1}{2}$ & 11.76\tabularnewline
\hline 
3 & 10 & 2 & 18.30\tabularnewline
\hline 
4 & 15 & $2\sfrac{1}{2}$ & 29.23\tabularnewline
\hline 
5 & 21 & 3 & 37.47\tabularnewline
\hline 
6 & 28 & $3\sfrac{1}{2}$ & 53.17\tabularnewline
\hline 
7 & 36 & 4 & 63.56\tabularnewline
\hline 
\end{tabular}
\par\end{center}%
\end{minipage}\hfill{} %
\begin{minipage}[t]{0.49\textwidth}%
\begin{center}
\begin{tabular}{c|c||c|c}
\hline 
Degree & \#Coeffs. & Ring & \#Nodes\tabularnewline
\hline 
\hline 
2 & 10 & 1 & 13.53\tabularnewline
\hline 
3 & 20 & $1\sfrac{1}{3}$ & 29.44\tabularnewline
\hline 
4 & 35 & $1\sfrac{2}{3}$ & 46.25\tabularnewline
\hline 
5 & 56 & 2 & 67.86\tabularnewline
\hline 
6 & 84 & $2\sfrac{1}{3}$ & 121.54\tabularnewline
\hline 
7 & 120 & $2\sfrac{2}{3}$ & 156.86\tabularnewline
\hline 
\end{tabular}
\par\end{center}%
\end{minipage}
\end{table}

\section{\label{sec:Overview-of-AES-FEM}Overview of AES-FEM}

Starting with a PDE with Dirichlet boundary conditions
\[
\begin{array}{cc}
\mathcal{L}u=f\quad & \text{on }\Omega\\
u=u_{D} & \text{on }\Gamma_{D}
\end{array}
\]
AES-FEM can be derived from (\ref{eq:varform-gwr-interior}) as follows. AES-FEM
uses generalized Lagrange polynomials as the basis functions $\left\{ \phi_{j}\right\} $
and the traditional FEM hat functions as the test functions $\left\{ \psi_{i}\right\} $.
The solution $u$ is approximated as $u=\sum u_{j}\phi_{j}$. As a
concrete example, consider the Poisson equation. For a given node
and corresponding test function $\psi_{i}$, we have
\begin{equation}
-\sum_{j=1}^{n}u_{j}\int_{\Omega}\nabla\psi_{i}\cdot\nabla\phi_{j}\ dV=\int_{\Omega}\psi_{i}f\ dV.\label{eq:aes-fem_weak_form}
\end{equation}

The stiffness matrix is assembled row by row. For each interior node,
a stencil is selected (see \ref{sec:Selection-of-Stencils}) and the
GLP basis functions are calculated on that stencil (see \ref{sec:Computation-of-GLP}).
If the stencil contains too few nodes, the stencil is enlarged (hence
the word \textit{adaptive} in the name adaptive extended stencil-FEM).
The adaptive expansion of the stencil ensures the stability of the
basis functions. The entries in the $i$th row of the stiffness matrix
and the $i$th entry of the load vector are calculated using (\ref{eq:aes-fem_weak_form}).
Dirichlet boundary conditions are enforced strongly.

The degree of the basis functions controls the order of convergence
of the method. Note that regardless of the degree of the basis functions,
AES-FEM always uses piecewise linear hat functions for the test functions,
so it requires only first-order meshes regardless of the degree of
its basis functions.

\section{\label{subsec:Functional-Analysis} Functional Analysis and Variational
Crimes}

Unlike FDM, of which the convergence follows from the fundamental
theorem of numerical analysis, proving convergence of FEM is more
complicated. It requires an intricate integration of functional analysis
and approximation theory, which were traditionally incompatible, and
hence the term ``variational crimes'' coined by Strang \cite{strang1972variational,strang1973analysis}.

\subsection{Functional Analysis of Coercive PDEs}

The convergence analysis of FEM is best known for coercive PDEs. One
of the most fundamental results is the \emph{Lax-Milgram lemma} \cite{lax1954parabolic}\cite[p. 83]{ern2013theory},
which states that an FEM is \emph{well-posed} (or \emph{invertible})
for \emph{bounded} and \emph{coercive} bilinear forms. Its proof boils
down to the \emph{Riesz representation theorem} \cite[p. 479]{ern2013theory}
for $C^{0}$ functions and the \emph{Poincaré inequality} \cite[p. 489]{ern2013theory}.
In practice, the boundedness and coercivity are satisfied on \emph{quasiuniform}
and \emph{well-shaped} meshes. In terms of convergence, for simpler
cases, such as FEM with Dirichlet boundary conditions, the error is
bounded in $L^{2}$ norm, for which the most successful technique
is the \emph{Aubin-Nitsche duality argument} \cite{aubin1967behavior,nitsche1968kriterium},
a.k.a. ``\emph{Nitsche's trick}'' \cite[p. 166]{strang1973analysis}.
 When approximation errors are involved, the convergence rates are
often proven only in $H^{1}$ norm (see e.g. \cite[p. 288]{brenner2008mathematical}
and \cite[p. 199]{ciarlet2002finite}).

\subsection{Functional Analysis of Noncoercive PDEs}

The generalization of functional analysis of FEM to noncoercive PDEs
requires the use of Banach or Sobolev spaces. The best known result
is the \emph{Banach-Nečas-Babuška} (\emph{BNB}) \emph{theorem} \cite[p. 84-85]{ern2013theory},
attributed to Nečas \cite{nevcas1962methode} and Babuška \cite{babuska1972survey},
regarding the \emph{invertibility} (or \emph{well-posedness}). It
generalizes the \emph{Lax-Milgram lemma}. The theorem states that
an FEM with a specific trial space $\Phi$ and test space $\Psi$
is \emph{invertible} if and only if 
\begin{equation}
\text{\ensuremath{\exists}}\alpha>0,\qquad\inf_{\phi\in\Phi\backslash\{0\}}\sup_{\psi\in\Psi\backslash\{0\}}\frac{a(\phi,\psi)}{\Vert\phi\Vert_{\Phi}\Vert\psi\Vert_{\Psi}}\geq\alpha\label{eq:inf-sup-condition}
\end{equation}
and 
\begin{equation}
\forall\psi\in\Phi,\qquad\sup_{\phi\in\Phi\backslash\{0\}}\vert a(\phi,\psi)\vert=0\Longrightarrow\psi=0,\label{eq:injective}
\end{equation}
where $\Vert\cdot\Vert_{\Phi}$ and $\Vert\cdot\Vert_{\Psi}$ are
some norms associated with the spaces $\Phi$ and $\Psi$ over $\Omega$,
respectively. An assumption of the BNB theorem is the \emph{boundedness}
of the bilinear form \cite[p. 82]{ern2013theory}
\begin{equation}
\text{\ensuremath{\exists}}C<\infty,\qquad\sup_{\phi\in\Phi\backslash\{0\}}\sup_{\psi\in\Psi\backslash\{0\}}\frac{a(\phi,\psi)}{\Vert\phi\Vert_{\Phi}\Vert\psi\Vert_{\Psi}}\leq C\label{eq:boundedeness}
\end{equation}
under some continuity requirements on $\Phi$ and $\Psi$ (such as
$C^{0}$ continuity). Eq.~(\ref{eq:inf-sup-condition}) is known
as the \emph{inf-sup condition}. For coercive problems, $\Vert\cdot\Vert_{\Phi}$
and $\Vert\cdot\Vert_{\Psi}$ in (\ref{eq:inf-sup-condition}) typically
correspond to some $L^{p}$ norm over $\Omega$. In practice, the
boundedness and inf-sup conditions also require \emph{quasiuniform}
and \emph{well-shaped} meshes. Since the invertibility condition is
purely algebraic, the solutions may suffer from spurious oscillations
for noncoercive PDEs.

\subsection{Variational Crimes\label{subsec:Variational-Crimes}}

In the classical functional analysis, a deviation from exact computations
or conforming FEM is considered a ``variational crime'' \cite{brenner2008mathematical,strang1973analysis}.
This includes interpolation errors that are not intrinsic in the $L^{p}$
(or $H^{1}$) norms, numerical integration errors, rounding errors,
etc. Among these, the most challenging is the geometric errors for
FEM with Neumann boundary conditions, which violates the assumptions
of Aubin-Nitsche duality argument. A more fundamental ``crime''
is the loss of continuity, which is introduced by \emph{nonconforming
finite elements} \cite[Section 10.3]{brenner2008mathematical}. This
loss of continuity violates the assumption of the Riesz representation
theorem, and its convergence analysis hence requires taking into account
the interface fluxes or jump conditions, even for coercive PDEs. AES-FEM
involves a similar ``crime'' due to its use of least-squares based
trial functions.
\end{document}